\documentclass[11pt]{article}
\usepackage{amsmath,amssymb,amsthm}
\usepackage[margin=1.1in]{geometry} 
\usepackage[colorlinks=true,urlcolor=blue]{hyperref}
\usepackage{graphicx}
\usepackage{booktabs}
\usepackage{listings}
\newtheorem{theorem}{Theorem}[section]
\newtheorem{proposition}[theorem]{Proposition}
\newtheorem{lemma}[theorem]{Lemma}
\newtheorem{corollary}[theorem]{Corollary}
\newtheorem{conjecture}[theorem]{Conjecture}
\theoremstyle{definition}
\newtheorem{definition}[theorem]{Definition}
\newtheorem{example}[theorem]{Example}
\newtheorem{remark}[theorem]{Remark}

\newcommand{\R}{\mathbb{R}}
\newcommand{\Sph}{\mathbb{S}}
\newcommand{\dd}{\,\mathrm{d}}
\newcommand{\Prob}{\mathcal{P}}
\newcommand{\Sgn}{\mathcal{M}_1}
\newcommand{\Mp}{M^{+}}
\newcommand{\Lell}{\mathcal{L}}
\DeclareMathOperator{\supp}{supp}
\DeclareMathOperator{\conv}{conv}

\title{When Positivity Matters:\\[2pt] Maximal Distance Integrals on Curves and Surfaces}
\author{Steven Finch\\
MIT Sloan School of Management}
\date{\today}

\begin{document}
\maketitle

\begin{abstract}
How far apart, on average, can $n$ points in a set $A$ be? As $n\to\infty$ the
answer is the \emph{maximal energy} $\Mp(A)$, the largest possible value of
$\mathbb{E}|X-Y|$ when $X$ and $Y$ are independent random points of $A$ with a
common distribution. If the distribution is allowed to carry negative mass, we
obtain a second constant $M(A)\ge \Mp(A)$. It is often far easier to compute: 
when a signed equilibrium measure exists, it is characterized by a 
\emph{linear} potential equation, whereas the
positive problem must also locate the support of the mass. This paper asks
when positivity matters, that is, when $\Mp(A)<M(A)$. For polygons and
polyhedra the extremal mass sits on the vertices and closed forms abound. Even
so, every hollow convex polygon has $\Mp<M$, and a simple kite shows how a
vertex can be forced out of the support. For smooth curves and surfaces the
answer depends on the shape. Circles and spheres have $\Mp=M$. The ellipse
$x^2+4y^2=4$ has $\Mp=M$ numerically, and we give its value to 50 digits.
The ellipses $x^2+9y^2=9$ and $x^2+25y^2=25$ have $\Mp<M$ numerically, with 
gaps in the support. For the ellipsoid $10x^2+15y^2+12z^2=15$,
studied by Larcher, Schmid and Wolf, we compute 50 digits of $M$ and give 
strong numerical evidence that $\Mp<M$. The extremal measure vanishes on 
two caps around the ends of the shortest axis. Proofs and numerical methods 
are collected in appendices.
\end{abstract}

\tableofcontents

\section{Introduction}\label{sec:intro}

\subsection{A puzzle about spreading points}

Scatter $n$ points in a compact set $A\subset\R^d$ so as to make the average
distance between pairs of them as large as possible. On the unit interval the
best strategy is plain: put half of the points at each end. The average
distance then tends to $\tfrac12$. What happens for a triangle, a square, a
circle, an ellipse? This paper is an invitation to that question.

\begin{definition}\label{def:gamma}
For compact $A\subset\R^d$, let
\[
\gamma_1(A)=\lim_{n\to\infty}\ \max_{\xi_1,\dots,\xi_n\in A}
\left(\frac{2}{n(n-1)}\sum_{j<k}|\xi_j-\xi_k|\right).
\]
\end{definition}

The limit exists and has a natural description in terms of mass
distributions. For a finite (signed) Borel measure $\sigma$ on $A$ define the
\emph{energy} and the \emph{potential}
\[
I(\sigma)=\iint|x-y|\dd\sigma(x)\dd\sigma(y),\qquad
U_\sigma(x)=\int|x-y|\dd\sigma(y).
\]
Let $\Prob(A)$ be the set of probability measures on $A$. If $X$ and $Y$ are
independent random points with law $\mu\in\Prob(A)$, then
$I(\mu)=\mathbb{E}|X-Y|$, the \emph{Gini mean difference} of $\mu$. Put
\[
\Mp(A)=\max_{\mu\in\Prob(A)} I(\mu).
\]

\begin{theorem}\label{thm:gamma}
For every compact $A$ and every $n\ge2$,
\[
\Mp(A)\ \le\ \max_{\xi_1,\dots,\xi_n\in A}\frac{2}{n(n-1)}\sum_{j<k}|\xi_j-\xi_k|
\ \le\ \frac{n}{n-1}\,\Mp(A).
\]
In particular $\gamma_1(A)=\Mp(A)$.
\end{theorem}

The idea is pleasantly probabilistic. For the left inequality, draw the points
independently from an optimal $\mu$: the \emph{expected} pair average is then
exactly $\Mp(A)$, so some configuration does at least that well. For the right
inequality, note that the $n$ points themselves define a probability measure.
Details are in Appendix~\ref{app:proofs}.

\subsection{Turning positivity off}

Now drop the requirement that the mass be nonnegative. Let $\Sgn(A)$ be the
set of finite signed Borel measures on $A$ of total mass $1$, and put
\[
M(A)=\sup_{\nu\in\Sgn(A)} I(\nu)\ \ \ge\ \Mp(A).
\]
Negative mass may look artificial. It makes the problem far more tractable,
because of one fact about the distance kernel: \emph{if $\sigma$ has total
mass zero, then $I(\sigma)\le0$, with equality only for $\sigma=0$}
(Lemma~\ref{lem:cnd}). Statisticians know this as the nonnegativity of the
energy distance~\cite{SR13}. Here are its consequences.

\begin{itemize}
\item Suppose a signed mass-one measure $\nu$ has \emph{constant} potential on
  $A$, $U_\nu\equiv\lambda$. Then $\nu$ is automatically the unique maximizer,
  and $M(A)=\lambda$ (Proposition~\ref{prop:signed}). We call $\nu$ the
  \emph{signed equilibrium measure}. Finding it is a linear problem; on a
  finite set it is the linear system $Dw=\lambda\mathbf 1$, where $D$ is the
  distance matrix. Such a measure need not exist: for the unit disk the supremum
  $M=\pi/2$ is not attained by any signed measure~\cite{HNW11}.
\item The positive problem adds the constraint $\mu\ge0$. This makes it an
  obstacle problem: one must find where the mass lives as well as how much
  there is.
\item The two problems agree exactly when the extremal probability measure
has constant potential, that is, when a signed equilibrium exists and is
a genuine probability measure (Theorem~\ref{thm:criterion}).
\end{itemize}

\noindent So the guiding question of this paper is:
\begin{center}
\emph{When does positivity matter, i.e.,\ when is $\Mp(A)<M(A)$?}
\end{center}
For the curves and surface we study, a signed equilibrium appears to exist,
with a density $\rho$ with respect to a natural parameter, and positivity
matters precisely when $\min\rho<0$.

\subsection{What we find}

\begin{itemize}
\item \emph{Polytopes} (Section~\ref{sec:poly}). The extremal probability
  measure lives on the vertices. This holds for the solid polytope, its
  boundary and its vertex set alike. Many closed forms follow: $\tfrac23$ for
  the equilateral triangle, $\tfrac{2+\sqrt2}{4}$ for the unit square,
  $\tfrac{3+3\sqrt2+\sqrt3}{8}$ for the unit cube, a formula for every
  triangle, and more.
\item \emph{The kite} (Section~\ref{sec:kite}). A slightly flattened kite
  shows how the signed equilibrium can give a vertex \emph{negative} mass. The
  true optimum then discards that vertex. A theorem of Wilson implies that the
  boundary of \emph{every} convex polygon has $\Mp<M$.
\item \emph{Circles and spheres} (Section~\ref{sec:circles}). By symmetry the
  uniform measure has constant potential, so $\Mp=M$, with values $4/\pi$ and
  $4/3$ for unit radius.
\item \emph{Ellipses} (Section~\ref{sec:ellipses}). For $x^2+4y^2=4$ the
  signed density is positive (numerically), so $\Mp=M$, which we compute to 50
  digits. For $x^2+9y^2=9$ and $x^2+25y^2=25$ the density changes sign. Our 
  computations show that $\Mp<M$ and find that the extremal measure abandons 
  the middle of the long flat sides altogether.
\item \emph{The ellipsoid} $10x^2+15y^2+12z^2=15$ (Section~\ref{sec:ellipsoid}).
  We compute 50 digits of $M$, give bounds for $\Mp$ that improve the
  estimate of Larcher, Schmid and Wolf~\cite{LSW00}, and find strong numerical
  evidence that $\Mp<M$.
\end{itemize}

\subsection{Why hollow sets?}\label{sec:hollow}

All the sets in this paper have empty interior: the kite is a polygonal
\emph{curve}, the ellipses are curves (1-manifolds), and the ellipsoid is a
surface (a 2-manifold). This is deliberate. For a convex body $K$ (compact and
convex, with interior points) two things happen automatically:
\begin{itemize}
\item the extremal probability measure lives on $\partial K$, so
  $\Mp(K)=\Mp(\partial K)$ (Theorem~\ref{thm:bjorck});
\item $\Mp(K)<M(K)$ (Corollary~\ref{cor:wilson}).
\end{itemize}
For example, the unit circle has $\Mp=M=4/\pi$, but the unit disk has
$\Mp=4/\pi$ and $M=\pi/2$, a supremum that no signed measure attains
\cite{HNW11}. Likewise the unit sphere has $\Mp=M=4/3$, while the unit ball
has $M=2$ \cite{Ale77,HNW11}.
So the solid tells us nothing new about $\Mp$, and $M$ is guaranteed to be
strictly larger. The interesting question lives on the hollow boundary.

\subsection{A word on names}\label{sec:names}

The constant $\Mp(A)$ has appeared under several names:
\begin{itemize}
\item a \emph{generalized transfinite diameter} of order one, in the tradition
  of Fekete and P\'olya--Szeg\H{o}~\cite{PS31};
\item Bj\"orck's \emph{maximal value of a generalized energy integral}~\cite{Bjorck56};
\item an \emph{optimal Riesz $(-1)$-energy}~\cite{BHS19};
\item a \emph{maximal mean distance} or maximal Gini mean difference.
\end{itemize}
By analogy with Newtonian capacity (kernel $|x-y|^{-1}$) it could also be called 
the \emph{Euclidean capacity} \cite{Finch19}. We simply say \emph{maximal energy}.

Our notation follows Hinrichs, Nickolas and Wolf~\cite{HNW11}, who wrote
$M(X)$ for the supremum over signed mass-one measures and $M^+(X)$ for the
supremum over probability measures. The notation elsewhere is not uniform:
Wolf~\cite{Wolf97} and Larcher--Schmid--Wolf~\cite{LSW00} wrote $M(X,d)$ for
the supremum over probability measures, our $\Mp$, while Nickolas and
Wolf~\cite{NW09,NW11a,NW11b,NWfin} used $M(X)$, as we do, for the signed
supremum. The signed problem goes back to Alexander and
Stolarsky~\cite{AS74}. A related but different constant, the rendezvous
number of Gross~\cite{Gross64}, is discussed in Appendix~\ref{app:rendezvous}.

\subsection{Summary of key numbers}\label{sec:summary}

In the table, ``ellipse'' means the curve and ``ellipsoid'' the surface;
$\rho$ is the density of the signed equilibrium measure, and ``conj'' marks a
numerically supported conjecture. Apart from the circle and the sphere, 
all entries are numerical. Values of $M$ come from spectral Galerkin computations; 
we print digits on which successive truncations agree, without a rigorous error bound. 
Bounds for $\Mp$ come from discrete optimization: the lower ends are energies of
explicit probability measures, and Appendix~\ref{app:rigor} explains how the
upper ends can be made rigorous. The entries ``$\min\rho<0$; $M^+<M$'' are
numerical conclusions in the same sense.

\begin{center}
	\renewcommand{\arraystretch}{1.25}
	\begin{tabular}{@{}p{4.3cm}p{8.6cm}p{2.6cm}@{}}
		\toprule
		Set & Value & Status\\
		\midrule
		Circle, radius $1$ & $M=4/\pi=M^{+}$ & proved\\
		Sphere, radius $1$ & $M=4/3=M^{+}$ & proved\\
Ellipse $x^2+4y^2=4$ &
\begin{tabular}{@{}l@{}} $M=2.108054212174101323143715462183\ldots$ \\
	$2.1080542101 < M^{+} < 2.1080542186$;\end{tabular} &
\begin{tabular}{@{}l@{}} $\min\rho\approx0.0073$; \\ $M^{+}=M$ conj \end{tabular}\\
Ellipse $x^2+9y^2=9$ & 
\begin{tabular}{@{}l@{}} $M=3.067727271436231046844229325280\ldots$ \\
	$3.0677272090 < M^{+} < 3.0677272273$ \end{tabular} &
\begin{tabular}{@{}l@{}} $\min\rho<0$; \\ $M^{+}<M$ \end{tabular}\\
Ellipse $x^2+25y^2=25$ & 
\begin{tabular}{@{}l@{}} $M=5.039416286675115620305281439677\ldots$ \\
	$5.0394159505 < M^{+} < 5.0394159819$ \end{tabular} &
\begin{tabular}{@{}l@{}} $\min\rho<0$; \\ $M^{+}<M$ \end{tabular}\\
\begin{tabular}{@{}l@{}} Ellipsoid \\ $10x^2+15y^2+12z^2=15$ \end{tabular} &
\begin{tabular}{@{}l@{}}$M=1.511269673044740064891740394794\ldots$ \\
	$1.5112600022 < M^{+} < 1.5112601265$ \\
	(LSW lower bound $1.5112509356$, $N=3969$) \end{tabular} &
\begin{tabular}{@{}l@{}}$\min\rho<0$; \\ $M^{+}<M$ \end{tabular} \\
		\bottomrule
	\end{tabular}
\end{center}

\section{Vertices win: polygons and polyhedra}\label{sec:poly}

\subsection{Where does the mass go?}

The basic structure theorem is due to Bj\"orck~\cite{Bjorck56}.

\begin{theorem}[Bj\"orck]\label{thm:bjorck}
Let $A\subset\R^d$ be compact.
\begin{enumerate}
\item There is a \emph{unique} $\mu_0\in\Prob(A)$ with $I(\mu_0)=\Mp(A)$.
\item Its potential satisfies $U_{\mu_0}\le\Mp(A)$ on $A$, with equality on
  $\supp\mu_0$.
\item If $A=K$ is a convex body and $d\ge2$, then $\supp\mu_0\subset\partial K$,
  and hence $\Mp(K)=\Mp(\partial K)$.
\end{enumerate}
\end{theorem}

Part~2 says that mass sits only where the potential is highest. Part~3 has a
simple reason: each function $x\mapsto|x-y|$ is convex, so every potential
$U_\mu$ is convex, and a convex function cannot have an interior maximum
unless it is constant. Convexity also tells us where the maximum lives on a
polytope, namely at a vertex:

\begin{theorem}[Polytopes reduce to vertices]\label{thm:polytope}
Let $K\subset\R^d$ be a convex polytope with vertex set $V$, and let $A$ be
compact with $V\subseteq A\subseteq K$. For example, $A$ may be $V$, the
boundary $\partial K$, or $K$ itself. Then $\Mp(A)=\Mp(V)$, and the extremal
measure of $A$ is that of $V$; in particular it lives on the vertices.
\end{theorem}

So for a polytope, its hollow boundary, or its vertex set, computing $\Mp$ is a
\emph{finite} problem: maximize $w^{\mathsf T}Dw$ over weight vectors $w\ge0$
with $\sum w_i=1$, where $D$ is the vertex distance matrix. Often symmetry
makes this immediate.

\begin{corollary}[Equal potentials]\label{cor:equalpot}
In the setting of Theorem~\ref{thm:polytope}, suppose weights $w_i\ge0$ with
$\sum_iw_i=1$ make the vertex potentials equal:
\[
U(v_i)=\sum_l w_l|v_i-v_l|=\lambda\quad\text{for every }i.
\]
Then $\Mp(A)=\lambda=\sum_{i,l}w_iw_l|v_i-v_l|$ and
$\mu_0=\sum_iw_i\delta_{v_i}$. Moreover, for every $n$ and all
$\xi_1,\dots,\xi_n\in A$,
\[
\sum_{j,k=1}^n|\xi_j-\xi_k|\ \le\ \lambda\,n^2,
\]
with equality when exactly $w_in$ of the points sit at each vertex $v_i$
(whenever these numbers are integers).
\end{corollary}

The last statement is a finite-$n$ fact, not merely a limit. For instance,
whenever $3\mid n$, putting $n/3$ points on each corner of an equilateral
triangle is \emph{exactly} optimal for the full double sum.

\begin{example}\label{ex:TSC}
Let $T$ be an equilateral triangle with unit sides, $S$ a unit square and $C$
a unit cube. Each may be taken solid, hollow, or as its vertex set. In each
case the symmetry group acts transitively on the vertices, so uniform weights
equalize the potentials, and $\lambda$ is the average distance from one vertex
to all vertices (including itself):
\begin{center}
\renewcommand{\arraystretch}{1.5}
\begin{tabular}{c l l l}
\toprule
Set & Distances from one vertex & $\Mp=\lambda$ & Value \\
\midrule
$T$ & $0,1,1$ & $\tfrac{2}{3}$ & $0.6666667$ \\
$S$ & $0,1,1,\sqrt2$ & $\tfrac{2+\sqrt2}{4}$ & $0.8535534$ \\
$C$ & $0,1,1,1,\sqrt2,\sqrt2,\sqrt2,\sqrt3$ & $\tfrac{3+3\sqrt2+\sqrt3}{8}$ & $1.1218364$ \\
\bottomrule
\end{tabular}
\end{center}
\end{example}

\subsection{A gallery of closed forms}\label{sec:gallery}

The same reasoning gives:

\begin{itemize}
\item \textbf{Rectangles and boxes.} The vertices are still transitive under
  the symmetries. For an $a\times b$ rectangle,
  $\Mp=\tfrac14\big(a+b+\sqrt{a^2+b^2}\big)$. For an $a\times b\times c$ box,
  \[
  \Mp=\tfrac18\Big(a+b+c+\sqrt{a^2+b^2}+\sqrt{b^2+c^2}
  +\sqrt{c^2+a^2}+\sqrt{a^2+b^2+c^2}\Big).
  \]
\item \textbf{Regular simplices.} With $m$ vertices and unit edges,
  $\Mp=(m-1)/m$. For the regular tetrahedron this is $\tfrac34$.
\item \textbf{Regular polygons.} For the regular $m$-gon with circumradius
  $1$, $\Mp=\frac{2}{m}\cot\frac{\pi}{2m}$. This equals $1.1547$, $1.2071$,
  $\dots$ for $m=3,4,\dots$, and increases to $4/\pi$, the value for the unit
  circle (Section~\ref{sec:circles}).
\item \textbf{Every triangle.} Let the sides be $a,b,c$, opposite the vertices
  $A,B,C$. The weights
  \[
  w_A\propto a(b+c-a),\qquad w_B\propto b(c+a-b),\qquad w_C\propto c(a+b-c)
  \]
  are nonnegative by the triangle inequality, and all three unnormalized
  vertex potentials equal $2abc$. Hence
  \[
  \Mp=\frac{2abc}{2(ab+bc+ca)-(a^2+b^2+c^2)} .
  \]
  This gives $\tfrac23$ for the equilateral triangle. In the degenerate case
  $a=1$, $b=c=\tfrac12$ it gives $\tfrac12$, the value for the unit interval.
\end{itemize}

\begin{remark}\label{rem:caution}
All of this concerns $\Mp$. For the vertex sets themselves the signed problem
has the same answer, since the weights above solve $Dw=\lambda\mathbf1$ and
are nonnegative. For the \emph{hollow} triangle $\partial T$, however,
$M(\partial T)>\Mp(\partial T)=\tfrac23$ (Corollary~\ref{cor:wilson}). The
triangle's boundary is our first example where positivity matters.
\end{remark}

\section{When a vertex must be left out: the kite}\label{sec:kite}

\subsection{Signed weights on finite sets}

For a finite set $V=\{v_1,\dots,v_m\}$ of distinct points with distance
matrix $D$, the signed problem is pure linear algebra. $D$ is invertible,
$\mathbf1^{\mathsf T}D^{-1}\mathbf1>0$, and
\[
M(V)=\frac{1}{\mathbf1^{\mathsf T}D^{-1}\mathbf1},\qquad
\text{attained by the weights } w=\frac{D^{-1}\mathbf1}{\mathbf1^{\mathsf T}D^{-1}\mathbf1}
\]
(Proposition~\ref{prop:finite}). These weights equalize the potential on $V$.
If they are all nonnegative, Corollary~\ref{cor:equalpot} applies and
$\Mp(V)=M(V)$. If one is negative, positivity matters.

Three points never produce a negative weight; that is the triangle formula
above. Rectangles, boxes and regular polygons cannot either, because symmetry
forces equal weights. So four points is the smallest case where a negative
weight can appear. The simplest example is a kite.

\subsection{The kite}

Take the kite with vertices
\[
A=(-1,0),\quad B=(1,0),\quad C=(0,1),\quad D=(0,-\varepsilon),
\qquad \varepsilon\ge0 \text{ small},
\]
and let $V=\{A,B,C,D\}$ and $K=\conv V$. The vertex $D$ sits just below the
midpoint of $AB$. Intuitively it is almost redundant: the long diagonal $AB$
already reaches past it on both sides.

\paragraph{The flat case $\varepsilon=0$.} The distances are $AB=2$,
$AC=BC=\sqrt2$ and $AD=BD=CD=1$. By symmetry $A$ and $B$ get the same weight.
Write the weights as $p,p,q,r$ at $A,B,C,D$, with $2p+q+r=1$. Then
\[
U(A)=2p+\sqrt2\,q+r,\qquad U(C)=2\sqrt2\,p+r,\qquad U(D)=2p+q .
\]
Setting $U(A)=U(D)$ gives $r=(1-\sqrt2)\,q$, which is negative whenever
$q>0$. Setting $U(A)=U(C)$ gives $q=(2-\sqrt2)\,p$, and normalizing yields
\[
p=\tfrac{2+\sqrt2}{8}\approx0.4268,\qquad q=\tfrac14,\qquad
r=\tfrac{1-\sqrt2}{4}\approx-0.1036,\qquad
\lambda=\tfrac{3+\sqrt2}{4}\approx1.1036 .
\]
The sign of $r$ is forced. Moving from $A$ to $D$ leaves the contribution of
$A$ and $B$ unchanged at $2p$, but $C$ is closer to $D$ than to $A$ ($1$
versus $\sqrt2$). The only way to keep the potential level is to give $D$
negative mass.

\paragraph{A genuine kite.} For $\varepsilon>0$ the quadrilateral is strictly
convex and $D$ is a true vertex. Now $AD=BD=\sqrt{1+\varepsilon^2}$ and
$CD=1+\varepsilon$. The weight at $D$ stays negative until it vanishes at the
root of
\[
2\sqrt{1+\varepsilon^2}+(2-\sqrt2)(1+\varepsilon)=2\sqrt2,
\qquad \varepsilon^*\approx0.2815 .
\]
A milder version of the same effect occurs for the three collinear points
$0,1,2$: solving $Dw=\mathbf 1$ gives the middle point weight exactly $0$.
Pulling one point off the line, as $C$ does here, turns that zero negative.

\subsection{What the optimum actually does}

For $0\le\varepsilon<\varepsilon^*$ the extremal probability measure drops
$D$. It is the measure of the right isosceles triangle $ABC$ (hypotenuse $2$,
legs $\sqrt2$), whose weights come from the triangle formula:
\[
w_A=w_B=\frac{\sqrt2}{4\sqrt2-2}\approx0.3867,\qquad
w_C=\frac{\sqrt2-1}{2\sqrt2-1}\approx0.2265,\qquad
\Mp=\frac{2}{2\sqrt2-1}\approx1.0938363 .
\]
To see that $D$ is rightly excluded, compute its potential. At
$\varepsilon=0.1$, for example,
$U(D)=2\sqrt{1.01}\,(0.3867)+1.1\,(0.2265)\approx1.026<1.0938$, so the
potential is at most $\Mp$ at all four vertices, as Theorem~\ref{thm:bjorck}
requires. At $\varepsilon=\varepsilon^*$, $U(D)$ rises to exactly $\Mp$; this
is the moment $D$ re-enters the support. By Theorem~\ref{thm:polytope} the
same value and measure serve for $V$, for the solid kite, and for the hollow
kite $\partial K$:
\[
\Mp(V)=\Mp(\partial K)=\Mp(K)=\frac{2}{2\sqrt2-1}\qquad(0\le\varepsilon<\varepsilon^*).
\]
One vertex of the polytope receives no mass at all.

\subsection{How much positivity costs}

The signed constant of the four-point set has a closed form. With
$s=\sqrt{1+\varepsilon^2}$ and $t=1+\varepsilon$,
\[
M(V)=\frac{t\,(2\sqrt2\,s-t)}{2\big[\sqrt2\,(s+t)+st-2-t^2\big]}
\]
(Appendix~\ref{app:kite}). A few values:
\begin{center}
\begin{tabular}{@{}cccc@{}}
\toprule
$\varepsilon$ & $M(V)$ & $\Mp(V)$ & signed weight at $D$\\
\midrule
$0$ & $\tfrac{3+\sqrt2}{4}=1.1035534$ & $1.0938363$ & $\tfrac{1-\sqrt2}{4}=-0.1036$\\
$0.1$ & $1.098586$ & $1.0938363$ & $-0.0706$\\
$\varepsilon^*\approx0.2815$ & $1.09384$ & $1.0938363$ & $0$\\
\bottomrule
\end{tabular}
\end{center}
At $\varepsilon^*$ the two constants meet, as they must once the signed
weight vanishes. For $\varepsilon\ge\varepsilon^*$ all four weights are
nonnegative and $\Mp(V)=M(V)$.

\subsection{Hollow polygons always pay}

Two general results complete the picture. The first characterizes when
positivity is free. It follows from Wolf's work~\cite{Wolf97}; a short proof
is in Appendix~\ref{app:signed}.

\begin{theorem}[Criterion]\label{thm:criterion}
For compact $A\subset\R^d$, $\Mp(A)=M(A)$ if and only if the extremal measure
$\mu_0$ has constant potential on all of $A$.
\end{theorem}

The second, due to David Wilson (see~\cite[Prop.~2]{CMY86} and~\cite{Wolf97}),
says that constant potentials are rare.

\begin{theorem}[Wilson]\label{thm:wilson}
Let $X\subset\R^d$ be compact and connected, and not a line segment. If some
probability measure on $X$ has constant potential on $X$, then no three points
of $X$ are collinear.
\end{theorem}

\begin{corollary}\label{cor:wilson}
If $X\subset\R^d$ is compact and connected, is not a segment, and contains
three collinear points, then $\Mp(X)<M(X)$. In particular:
\begin{enumerate}
\item the boundary of every convex polygon, e.g.,\ the hollow triangle
  $\partial T$, the hollow square $\partial S$ and the hollow kite $\partial K$;
\item every convex body in $\R^d$, $d\ge2$.
\end{enumerate}
\end{corollary}

For polygonal curves there is a short direct proof
(Lemma~\ref{lem:polygon}): along an edge, a positive measure's potential bends
too much to stay level. No computation is needed. For the kite the gap can
also be measured. Since $V\subset\partial K$, every signed measure on $V$ is a
signed measure on $\partial K$, so
\[
M(\partial K)-\Mp(\partial K)\ \ge\ M(V)-\frac{2}{2\sqrt2-1}\ >\ 0
\qquad(0\le\varepsilon<\varepsilon^*),
\]
about $0.0097$ at $\varepsilon=0$ and $0.0048$ at $\varepsilon=0.1$. We have 
not computed $M(\partial K)$ itself, but it is finite, thanks to a
simple bound obtained by projecting onto lines, as Hinrichs, Nickolas and
Wolf do for the ball~\cite{HNW11}.

\begin{proposition}[Width bound]\label{prop:width}
	Let $X\subset\R^2$ be compact and let $P$ be the perimeter of its convex
	hull. Then $M(X)\le P/4$. More generally, for $X\subset\R^n$,
	$M(X)\le\tfrac12\,b(n)\,\bar w$, where $\bar w$ is the mean width of
	$\conv X$ and $b(n)=\sqrt\pi\,\Gamma(\frac{n+1}2)/\Gamma(\frac n2)$.
	Equality holds for segments and for solid balls~\cite{HNW11}.
\end{proposition}

\noindent The proof is in Appendix~\ref{app:width}. For example, the hollow
unit square satisfies $\tfrac{2+\sqrt2}{4}<M(\partial S)\le1$, and the hollow
unit triangle $\tfrac23<M(\partial T)\le\tfrac34$. For the kite,
$M(\partial K)\le\tfrac12\big(\sqrt2+\sqrt{1+\varepsilon^2}\big)$, which is
about $1.207$ at $\varepsilon=0$. The bound is far from sharp for smooth
curves: the unit circle has $M=4/\pi\approx1.273$ against $P/4=\pi/2$.

Polygons, then, always pay a price for positivity, and the reason is
geometric: they have straight edges. Smooth curves have no collinear triples,
so Wilson's theorem says nothing about them. The rest of the paper explores
what happens there.

\section{Circles and spheres}\label{sec:circles}

On a circle or a sphere every point looks like every other, so the uniform
measure has constant potential. By Theorem~\ref{thm:criterion} and
Proposition~\ref{prop:signed}, the uniform measure is then extremal for both
problems:
\[
\Mp=M=\frac{1}{2\pi}\int_0^{2\pi}2R\left|\sin\tfrac\theta2\right|\dd\theta=\frac{4R}{\pi}
\quad\text{(circle of radius $R$)},
\]
\[
\Mp=M=\frac R2\int_{-1}^1\sqrt{2-2t}\dd t=\frac{4R}{3}
\quad\text{(sphere of radius $R$)}.
\]
For the sphere, $t=\langle u,v\rangle$ is uniformly distributed on $[-1,1]$
when $v$ is uniform. These classical facts go back to Bj\"orck~\cite{Bjorck56}
and Stolarsky~\cite{Sto73}; see also Morris--Nickolas~\cite{MN83} for spheres
in all dimensions. The filled disk has the same $\Mp=4/\pi$: the potential of
the uniform measure on the circle is convex, so it is at most $4/\pi$ on the
disk, and Proposition~\ref{prop:bracket} below gives the upper bound.

The regular $m$-gons of Section~\ref{sec:gallery} approach the circle from
below, $\frac2m\cot\frac\pi{2m}\uparrow\frac4\pi$. Their vertex sets have
$\Mp=M$, while their hollow boundaries have $\Mp<M$ by
Corollary~\ref{cor:wilson}.

\section{Ellipses}\label{sec:ellipses}

\subsection{Setting the stage}

Write $E_{a,b}$ for the ellipse \emph{curve} $x(t)=(a\cos t,\,b\sin t)$,
$0\le t<2\pi$, with $a>b$. We regard $t$ as an angle, i.e.,\ modulo $2\pi$; 
accordingly $|t|$ and $|t-\pi|$ denote distances on the circle 
$\R/2\pi\mathbb{Z}$, so that, for example, 
$\{|t|\le\tau\}=[0,\tau]\cup[2\pi-\tau,2\pi)$.
Thus $E_{2,1}$ is $x^2+4y^2=4$, $E_{3,1}$ is
$x^2+9y^2=9$, and $E_{5,1}$ is $x^2+25y^2=25$. Scaling multiplies both
constants: $\Mp(cA)=c\,\Mp(A)$ and $M(cA)=c\,M(A)$. So only the ratio $b/a$
matters. By Theorem~\ref{thm:bjorck} the filled ellipse has the same $\Mp$ as
the curve. It has $\Mp<M$ automatically (Corollary~\ref{cor:wilson}), which is
why we study the curve.

Ellipses are genuinely harder than everything so far:
\begin{itemize}
\item there are no vertices to hold the mass, and the symmetry group is
  finite, so no symmetry dictates the answer;
\item the distance kernel factors as
  \[
  |x(t)-x(s)|=2\left|\sin\tfrac{t-s}{2}\right|
  \sqrt{a^2\sin^2\tfrac{t+s}{2}+b^2\cos^2\tfrac{t+s}{2}},
  \]
  a function of $t-s$ times a function of $t+s$, so Fourier modes do not
  diagonalize it;
\item we know of no closed form for the answer, and do not expect one.
\end{itemize}
We describe densities with respect to the parameter $t$. Their mean value is
$1/(2\pi)\approx0.159$.

\subsection{Brackets}

Numerical answers become trustworthy through one simple inequality, proved in
Appendix~\ref{app:signed}.

\begin{proposition}[Bracketing]\label{prop:bracket}
For every signed mass-one measure $\nu$ on $A$,
\[
\Mp(A)\ \le\ 2\max_A U_\nu-I(\nu).
\]
If $\nu$ is a probability measure, then also $I(\nu)\le\Mp(A)$.
\end{proposition}

So every trial probability measure traps $\Mp$ between two computable
numbers. Here is a first attempt on $E_{2,1}$:
\begin{center}
\begin{tabular}{@{}lc@{}}
\toprule
Probability measure on $E_{2,1}$ & Energy\\
\midrule
half the mass at each end of the major axis & $2$\\
uniform in $t$ & $2P/\pi^2\approx1.96330$\\
best measure on the four axis endpoints & $2.037854$\\
the extremal measure $\mu_0$ & $2.10805421217\ldots$\\
\bottomrule
\end{tabular}
\end{center}
Here $P\approx9.68845$ is the perimeter. Uniform mass does worse than two
points! For the four-point measure, put mass $p$ at $(\pm2,0)$ and
$q=\tfrac12-p$ at $(0,\pm1)$. Then $I=8p^2+4q^2+8\sqrt5\,pq$, which is
maximized at $p^*\approx0.41982$. Its potential at the parameter values
$t=0,\ \pi/6,\ 0.6,\ \pi/4,\ \pi/3$ is about
$2.0379,\ 2.1472,\ 2.1488,\ 2.1428,\ 2.1170$. The potential exceeds the
energy away from the vertices, so by Theorem~\ref{thm:bjorck} the four-point
measure is not optimal. Proposition~\ref{prop:bracket} gives roughly
$\Mp\lesssim2(2.1488)-2.0379\approx2.26$. Much finer trial measures, and
spectral methods, narrow this bracket to the values in
Section~\ref{sec:summary}.

\subsection{The ellipse $x^2+4y^2=4$: positivity is free}\label{sec:E21}

This is the ellipse treated by Larcher, Schmid and Wolf~\cite[p.~313]{LSW00}.
A spectral Galerkin computation (Appendix~\ref{app:spectral-ellipse}) gives 
the signed constant; its first 50 digits, on which successive truncations 
agree, are:
\[
M(E_{2,1})=2.1080542121741013231437154621832014783797802837239\ldots
\]
The signed density has its maximum, about $0.52$, at the ends of the major
axis ($t=0,\pi$). It dips to its minimum
\[
\min_t\rho(t)=0.0073289961\ldots>0
\]
at the ends of the minor axis ($t=\pm\pi/2$), about $1/70$ of the peak. Since
the signed equilibrium is positive, it \emph{is} the extremal probability
measure. Therefore $\Mp(E_{2,1})=M(E_{2,1})$ and $\mu_0$ has full support.
Independent dense computations give the bracket in the summary table, which
agrees with $M$ to eight digits.

This conclusion rests on numerics. A rigorous proof would need a bound on the
truncation error of $\rho$ smaller than $0.0073$. The margin is huge compared
with the observed errors, below $10^{-20}$, so a computer-assisted proof looks
well within reach (Section~\ref{sec:questions}).

\paragraph{Comparison with LSW.} Larcher, Schmid and Wolf obtained
$2.108054066642\ldots$ from a 4095-point grid uniform in $t$. Our dense
computations show that the optimum on such a grid lies below $\Mp$ by about
$2.4404/N^2$. For $N=4095$ this predicts $2.108054066650$, in agreement with
their value to about $10^{-11}$. So their computation is accurate for the grid
problem, but only its first seven significant digits are digits of $\Mp$.
Their density histogram, peaking at about $0.52$ at $t=0,\pi$, matches ours.
Their conjecture that the extremal measure has full support is supported by
our computation.

\subsection{The ellipse $x^2+25y^2=25$: positivity matters}\label{sec:E51}

Flatten the ellipse to $b/a=\tfrac15$ and the picture changes.

\paragraph{The signed equilibrium dips below zero.} The spectral computation
converges quickly (Appendix~\ref{app:spectral-ellipse}) and gives
\[
M(E_{5,1})=5.0394162866751156203052814396774901490677251309166\ldots
\]
The signed density rises to tall peaks of about $1.42$ at $t=0,\pi$. Along
the long flat sides it is almost zero, and it is slightly \emph{negative} on
one band in each half of the curve:
\[
\rho<0\quad\text{on}\quad (t_-,\pi-t_-)\cup(\pi+t_-,2\pi-t_-),\qquad t_-=0.779018 .
\]
Each band contains an end of the minor axis and corresponds to $|x|/a<0.7116$.
Within a band the density has a shallow ``W'' shape. Its two minima,
$\rho=-0.000411098$, are at $t=0.911007$ and $\pi-0.911007$, i.e.,\ at
$x/a=\pm0.612$. In between it has a local maximum,
$\rho(\pi/2)=-0.0000797838$, at the end of the minor axis. The negative part
is tiny: $0.26\%$ of the mean density and $0.03\%$ of the peak. Still, it 
exceeds the apparent truncation error by many orders of magnitude,
and if genuine it means that positivity matters (Corollary~\ref{cor:gap}).

\paragraph{The extremal measure abandons the flat sides.} Dense 
computations with up to $30000$ nodes (Appendix~\ref{app:dense}) give
\[
5.0394159505\ <\ \Mp(E_{5,1})\ <\ 5.0394159819,\qquad
\Mp(E_{5,1})\approx5.0394159581 .
\]
The upper bound lies $3.05\times10^{-7}$ below $M$, so $\Mp<M$ is
established numerically, independently of the sign of $\rho$. 
Appendix~\ref{app:rigor} explains how this comparison can be made 
fully rigorous at modest cost.The estimated difference is 
$M-\Mp\approx3.3\times10^{-7}$.

The support of $\mu_0$ consists of two arcs about the ends of the major axis:
\[
\supp\mu_0=\{\,|t|\le t_e\,\}\cup\{\,|t-\pi|\le t_e\,\}
=[0,t_e]\cup[\pi-t_e,\pi+t_e]\cup[2\pi-t_e,2\pi),\qquad t_e\approx0.632,
\]
that is, the parts of the curve with $|x|/a\ge0.807$. These arcs account for
$40\%$ of the parameter range but only about $23\%$ of the length. The two
flat sides carry no mass at all. On the support the density of $\mu_0$ closely 
matches the signed density, right up to the edge, where it is still about 
$0.004$ (a quarter of a percent of its peak); there it drops abruptly to zero. 
Numerically, then, the density of $\mu_0$ has a small jump at the free boundary, 
with no spike or rim (Figure~\ref{fig:rho-zoom}).

\paragraph{The gap is wider than the negative region.}
\begin{center}
\begin{tabular}{@{}lll@{}}
\toprule
 & where $\rho<0$ & gap, where $\mu_0$ has no mass\\
\midrule
$t$-interval (first half) & $(0.779,\,2.363)$ & $(0.632,\,2.509)$\\
fraction of the $t$-range & $0.504$ & $0.597$\\
$|x|/a$ & $<0.712$ & $<0.807$\\
\bottomrule
\end{tabular}
\end{center}
On the arcs $0.632<|t|<0.779$ (and their mirror images) the signed density is
\emph{positive}, yet $\mu_0$ puts no mass there. The support of $\mu_0$ is not
simply ``where $\nu$ is positive''; it is decided by the potential. Removing
the negative mass lowers the potential below $\Mp$ on a wider set than the
one where $\nu<0$. For Riesz $s$-kernels with $d-2\le s<d$, Brauchart,
Dragnev and Saff~\cite{BDS09} proved that the support of the extremal measure
is contained in that of the positive part of the signed equilibrium. We know
of no such theorem for the kernel $|x-y|$, though both of our examples obey
the inclusion.

\subsection{The ellipse $x^2+9y^2=9$ and the threshold}\label{sec:E31}

The intermediate ellipse $E_{3,1}$ behaves like $E_{5,1}$, but less
extremely. The spectral computation gives
\[
M(E_{3,1})=3.0677272714362310468442293252800224326858727568806\ldots,
\]
the truncations $K=80$ and $K=160$ agreeing in all 50 digits. The signed density again 
numerically dips below zero, with $\min\rho\approx-0.0005105$, hence positivity matters. 
Unlike $E_{5,1}$, the signed density here decreases monotonically along the
flat side: it is negative on a single band $t_-<t<\pi-t_-$, $t_-\approx1.30$,
about each end of the minor axis, and attains its minimum at $t=\pi/2$. This
band again lies inside the gap of $\supp\mu_0$, and on $1.124<t<1.30$ the
signed density is positive although $\mu_0$ carries no mass. Dense computations give
\[
3.0677272090<\Mp(E_{3,1})<3.0677272273,\qquad \Mp(E_{3,1})\approx3.0677272134,
\]
so $M-\Mp\approx5.8\times10^{-8}$. The upper bound at $n=30000$ lies $4.4\times10^{-8}$ 
below $M$, so $\Mp<M$ is established numerically, independently of the sign of $\rho$.
The extremal density peaks at about $0.83$
at the ends of the major axis. It is supported on two arcs $|t|\le t_e$ and
$|t-\pi|\le t_e$ with $t$ modulo $2\pi$ and $t_e\approx1.124$, 
i.e.,\ where $|x|/a\ge0.43$. 
These arcs hold $72\%$ of the parameter range and about $61\%$ of the length.
Compared with $E_{5,1}$, the gaps are narrower and confined to the middle of
the flat sides.

Curiously, the minimum of $\rho$ is \emph{not} monotone in $b/a$:
\begin{center}
	\begin{tabular}{@{}ccccccc@{}}
		\toprule
		$b/a$ & $\tfrac12$ & $\tfrac25$ & $\tfrac13$ & $\tfrac14$ & $\tfrac15$ & $\tfrac1{10}$\\
		\midrule
		$\min\rho$ & $+0.0073$ & $+0.00039$ & $-0.00051$ & $-0.00036$ & $-0.00041$ & $-0.00078$\\
		\bottomrule
	\end{tabular}
\end{center}
The minimum falls through zero between $b/a=\tfrac25$ and $\tfrac13$, rises
again to a local maximum near $b/a=\tfrac14$, and then deepens. A plausible
explanation is that two minima compete: one at the end of the minor axis,
which dominates for $b/a\ge\tfrac13$, and the off-axis minima of the ``W''
shape, which dominate by $b/a=\tfrac15$. For thinner ellipses the density
concentrates in ever taller peaks at the ends of the major axis (of height
roughly proportional to $a/b$), and the negative dip appears to grow with
them. The depth of the dip is in any case a poor guide to the size of
$M-\Mp$: $E_{3,1}$ has a deeper dip than $E_{5,1}$, yet its energy gap is
about six times smaller.

\bigskip\bigskip

\begin{figure}[ht]
	\centering
	\includegraphics[width=0.8\textwidth]{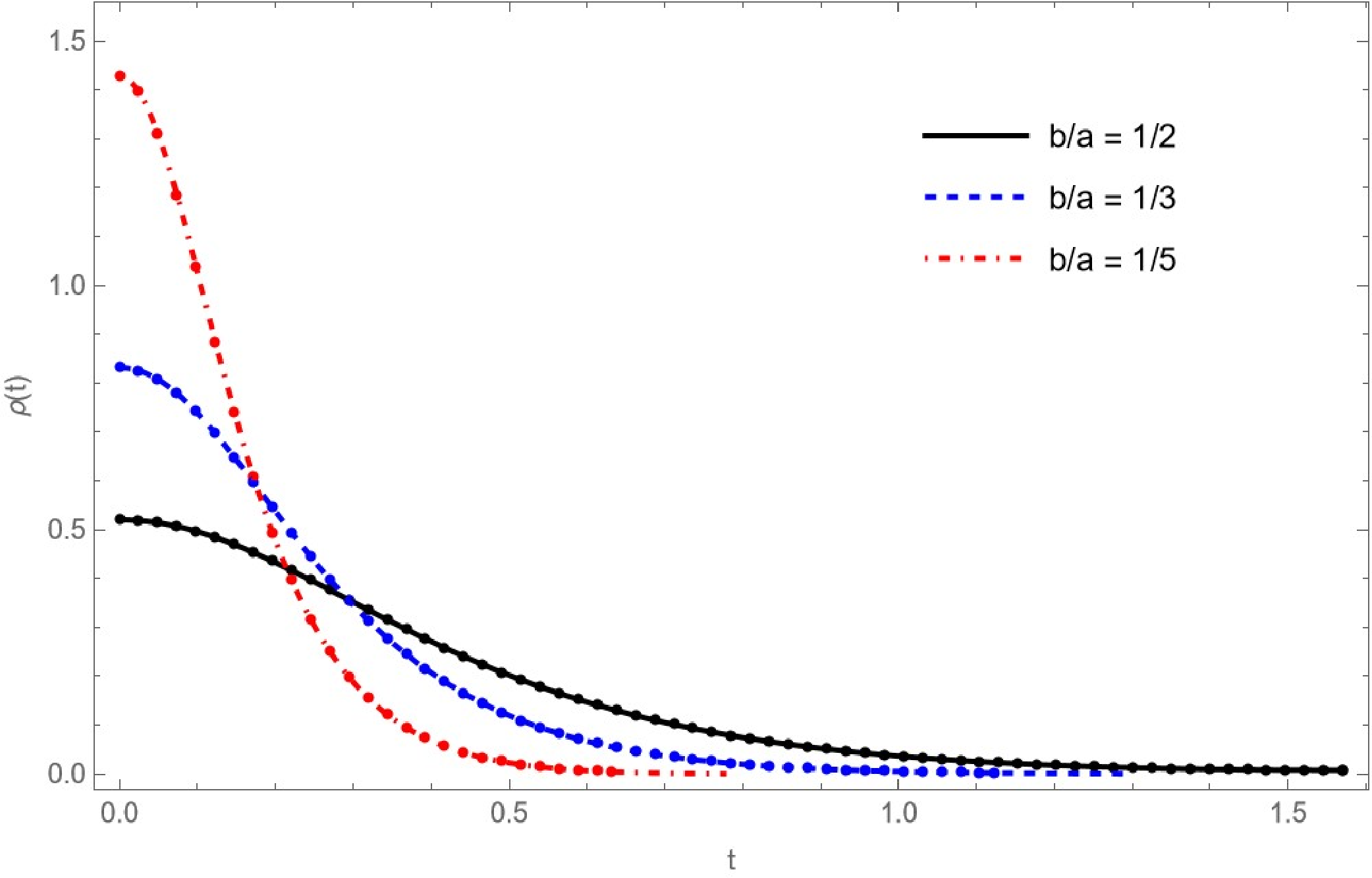}
	\caption{Densities with respect to $t$ on the quarter-curve $0\le t\le\pi/2$
		for $b/a=\tfrac12,\tfrac13,\tfrac15$. Curves: the signed equilibrium density
		$\rho$ (spectral, $K=80$). Dots: the extremal probability measure $\mu_0$
		(dense computation, $n=6400$). For $b/a=\tfrac12$ the two coincide and
		$\mu_0$ has full support; for $b/a=\tfrac13$ and $\tfrac15$ the dots stop at
		the edge of the support, $t_e\approx1.124$ and $0.632$.}
	\label{fig:rho-shape}
\end{figure}

\clearpage

\begin{figure}[ht]
	\centering
	\includegraphics[width=0.8\textwidth]{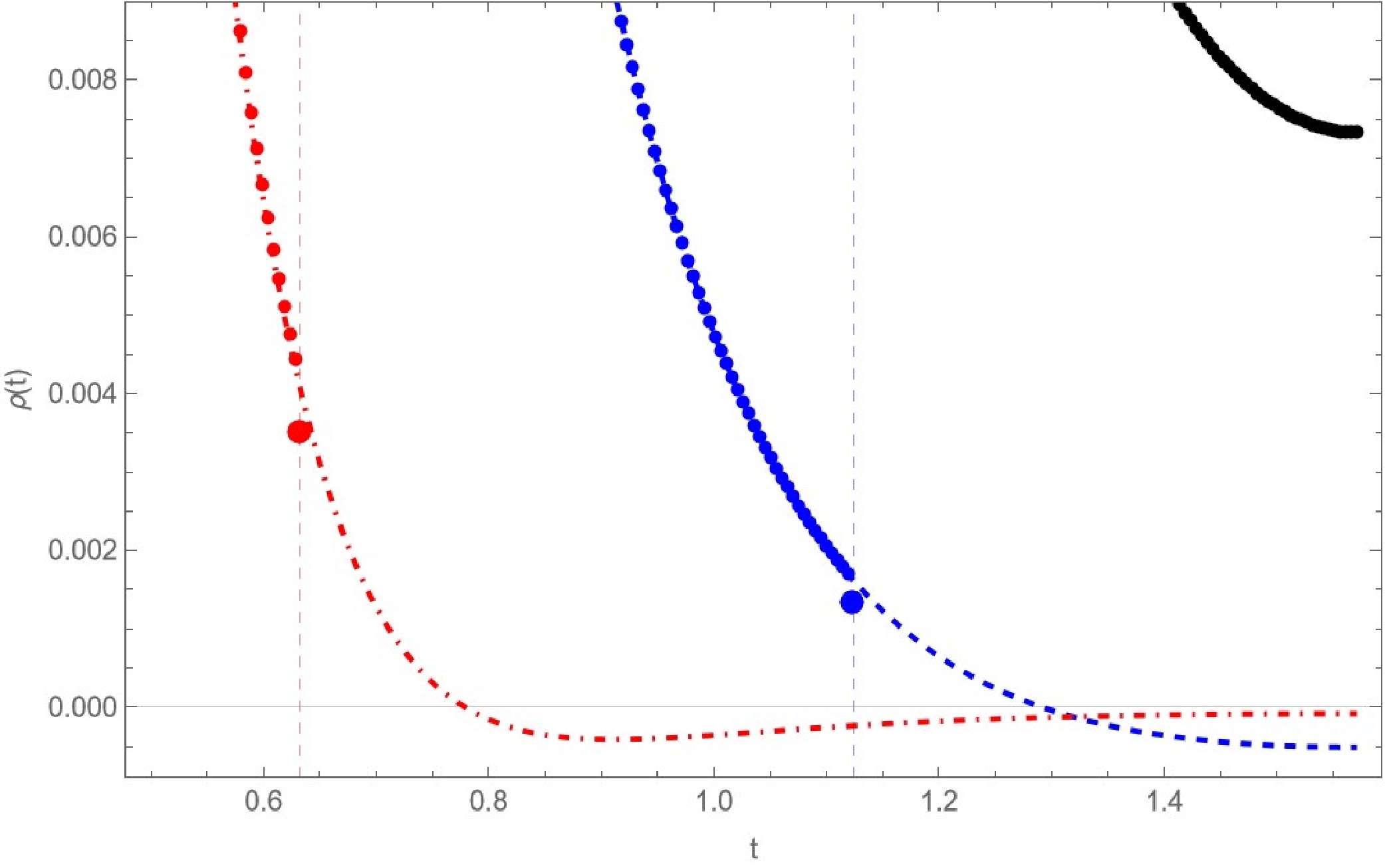}
	\caption{The same densities near zero. For $b/a=\tfrac12$ the signed density
		stays positive (minimum $0.0073$ at $t=\pi/2$). For $b/a=\tfrac13$ and
		$\tfrac15$ it dips below zero; the dip is deeper for $\tfrac13$, although the
		energy gap $M-\Mp$ is larger for $\tfrac15$. Dashed lines mark the support
		edges; between each edge and the nearest sign change, $\rho>0$ yet $\mu_0$
		carries no mass. For $b/a=\tfrac13$ the minimum is at $t=\pi/2$; 
		for $b/a=\tfrac15$ the dip has a W shape, with minima at $t=0.911$. 
		At the edge of the support the density of $\mu_0$ (dots) drops from 
		a small positive value to zero.}
	\label{fig:rho-zoom}
\end{figure}

\

\subsection{From ellipse to segment}

As $b/a$ decreases, the extremal measure passes through three stages:
\begin{itemize}
\item at $b/a=\tfrac12$ it has full support and a smooth positive density,
  bunched toward the ends of the major axis;
\item at $b/a=\tfrac13$ it lives on two arcs covering about $61\%$ of the
  length (peak density $0.83$), and at $b/a=\tfrac15$ on two shorter arcs
  covering about $23\%$ (peak $1.42$); in both cases its density falls to a 
  small positive value (about $0.0017$ and $0.004$, respectively) and then 
  drops to zero at the edges;
\item as $b/a\to0$ it tends to $\tfrac12(\delta_{(-a,0)}+\delta_{(a,0)})$, the
  extremal measure of a segment, and $\Mp\to a$.
\end{itemize}
Already at $b/a=\tfrac15$ the ellipse behaves like a ``thick segment'':
$\Mp(E_{5,1})/5=1.0078832$, against $1$ for the segment. Somewhere between 
$b/a=\tfrac13$ and $\tfrac25$ lies a critical ratio $\kappa^*$ at which 
the signed density first touches zero.

\section{The ellipsoid $10x^2+15y^2+12z^2=15$}\label{sec:ellipsoid}

\subsection{The surface}

Following Larcher, Schmid and Wolf~\cite[Example~3]{LSW00}, let $\Lell$ be the
ellipsoid \emph{surface}
\[
10x^2+15y^2+12z^2=15,\qquad\text{i.e.,}\qquad
\frac{x^2}{3/2}+y^2+\frac{z^2}{5/4}=1,
\]
with semi-axes $(a,b,c)=(\sqrt{3/2},\,1,\,\sqrt5/2)\approx(1.2247,\,1,\,1.1180)$.
It is a gently squashed sphere: the axes differ by at most $22\%$. We
parametrize it by the unit sphere, $x=(au_1,bu_2,cu_3)$ with $u\in\Sph^2$, and
write $\rho(u)$ for densities with respect to the area measure of
$\Sph^2$. The uniform measure in $u$ has energy $1.48871$, a first lower bound
for $\Mp$.

\subsection{The signed equilibrium: 50 digits and a sign change}

A spectral Galerkin method with spherical harmonics
(Appendix~\ref{app:spectral-ellipsoid}) converges by about one digit per
harmonic degree. Degrees $48$ and $56$ agree to $51$ significant digits, and
we report
\[
M(\Lell)=1.5112696730447400648917403947947430276305117569802\ldots
\]
These digits are stable under changes of degree, quadrature and working
precision, but we have no rigorous error bound yet.

Numerically, even this mild ellipsoid has a signed density that changes sign. The 
density is smallest exactly at $(0,\pm1,0)$, the ends of the \emph{shortest} axis,
where $\rho_{\min}=-0.0169305$. (The density with respect to surface area is
$-0.0123643$ there.) The region $\{\rho<0\}$ consists of two slightly
non-circular caps about these points. Their largest angular radius is
$27.4^\circ$, and they cover $10.6\%$ of the parameter sphere and $11.5\%$ of
the surface area. Compared with the ellipses, the negative part is large: its
depth is about a fifth of the mean density $1/(4\pi)$. A heuristic reason is
that inverting the distance operator amplifies the quadrupole ($l=2$)
component of the shape by a factor of $35$.

\subsection{The extremal probability measure}

Discrete optimization on octant-symmetric point sets with up to
$240{,}000$ effective points (Appendix~\ref{app:qp}) gives
\[
1.5112600022\ <\ \Mp(\Lell)\ <\ 1.5112601265,\qquad
\Mp(\Lell)\approx1.511260014 .
\]
The lower bound is the energy of an explicit probability measure. The upper
bound comes from a grid search followed by local maximization of the
potential, and is not rigorous (Appendix~\ref{app:rigor}). An error in the
search, however, would change it by roughly the bracket width
($\sim10^{-7}$), far less than the computed separation
\[
M(\Lell)-\Mp(\Lell)\ \ge\ 9.55\times10^{-6}
\]
(estimated value $9.66\times10^{-6}$). Independently, the computed signed
density is negative by a wide margin, which by Corollary~\ref{cor:gap} would
imply $\Mp<M$; but its sign is also a numerical result. Together these give
strong numerical evidence that positivity matters. The extremal
measure vanishes on two caps about $(0,\pm1,0)$ of angular radius about
$41.5^\circ$, covering about a quarter of the parameter sphere. As for the
ellipse $E_{5,1}$, these gaps are markedly larger than the negative caps
($27^\circ$) of the signed density. The signed solution locates the gaps but
underestimates their size.

\subsection{Comparison with Larcher--Schmid--Wolf}

LSW maximized over probability measures~\cite[pp.~296--297]{LSW00}. With
$N=3969$ grid points they reported $1.51125093560193590$, a lower bound for
$\Mp(\Lell)$. Their nodes form a $63\times63$ latitude--longitude grid in
which all $63$ points of the first row coincide at the north pole, so only
$3907$ points are distinct. Their own reported potentials show that their
value is within $2\times10^{-8}$ of the optimum on that grid. The shortfall of
about $9.1\times10^{-6}$ relative to $\Mp$ therefore comes from the placement
of the nodes, not from slow convergence of their algorithm.

Their data independently confirm the gaps. On the ellipse $E_{2,1}$ their node
potentials spread over only $3.7\times10^{-8}$. On the ellipsoid they spread
over $1.36\times10^{-3}$, far too much for the potential to be constant on a
full support. Their Fig.~3 shows two regions in which their algorithm places
essentially no mass, consistent with the two gaps found here.

\section{Questions for further exploration}\label{sec:questions}

\begin{enumerate}
\item \textbf{The critical ellipse.} Find the ratio $\kappa^*\in(\frac13,\frac25)$
  at which the signed density of the ellipse first touches zero. The minimum of 
  $\rho$ lies at $t=\pi/2$ for $b/a=\tfrac12$ and $\tfrac13$, and off the axis 
  only by $b/a=\tfrac15$, which suggests the first zero appears at the ends of 
  the minor axis. Does the
  first zero appear at the ends of the minor axis, as the shape at
  $b/a=\tfrac12$ suggests, or off the axis, as at $b/a=\tfrac15$?
  \begin{conjecture}
  There is $\kappa^*\in(\tfrac13,\tfrac25)$ such that $\Mp(E_{a,b})=M(E_{a,b})$
  for $b/a\ge\kappa^*$ and $\Mp(E_{a,b})<M(E_{a,b})$ for $b/a<\kappa^*$.
  \end{conjecture}
\item \textbf{Thin ellipses.} How does $\min\rho$ behave as $b/a\to0$? Does
the signed density near the ends of the major axis approach a fixed profile
after rescaling $t$ by $a/b$?
\item \textbf{A proof for $E_{2,1}$.} Certify $\rho>0$ rigorously, for example
  by an interval-arithmetic bound $\|\rho-\rho_K\|_\infty<0.0073$ on the
  truncated Fourier system.
\item \textbf{Ellipsoids.} For a family such as semi-axes $(1,1,c)$, find where
  the signed density first changes sign. Below that threshold, $\Mp=M$ is
  computable to many digits.
\item \textbf{Supports and negative regions.} Is the region where the signed
  density is negative always contained in the gap of $\supp\mu_0$, for the
  kernel $|x-y|$?
\item \textbf{High precision for $\Mp$.} Develop a free-boundary spectral
  method that yields $\Mp$ to many digits when $\Mp<M$.
\item \textbf{Edges of the support.} Numerically the density of $\mu_0$ jumps
from a small positive value to zero at the edge of its support on $E_{3,1}$
and $E_{5,1}$. Prove this, and determine the behavior on $\Lell$.
\item \textbf{Polygons and convex bodies.} By Proposition~\ref{prop:width},
$\tfrac{2+\sqrt2}{4}<M(\partial S)\le1$ for the hollow unit square. What is
$M(\partial S)$, and is it attained? More boldly: does every planar convex
body $K$ satisfy $M(K)=P/4$, as the disk does?
\item \textbf{Polygonal curves.} Compute $M(\partial K)$ for a convex polygon.
  Is it finite?
\item \textbf{Closed forms.} Is $\Mp(E_{2,1})$ expressible through classical
  constants, such as complete elliptic integrals?
\end{enumerate}

\appendix

\section{Proofs}\label{app:proofs}

Throughout, $A\subset\R^d$ is compact, and measures are finite signed Borel
measures on $A$.

\subsection{Negative definiteness of the distance kernel}

\begin{lemma}[Crofton]\label{lem:crofton}
For $x\in\R^d$ and the (unnormalized) surface measure $\sigma$ on
$\Sph^{d-1}$,
\[
\int_{\Sph^{d-1}}|\langle x,\omega\rangle|\dd\sigma(\omega)=c_d\,|x|,
\qquad c_2=4,\quad c_3=2\pi .
\]
\end{lemma}
\begin{proof}
By rotation invariance we may take $x=|x|e_d$. For $d=2$,
$\int_0^{2\pi}|\cos\theta|\dd\theta=4$. For $d=3$, with
$t=\langle\omega,e_3\rangle$ we have $\dd\sigma=\dd t\dd\phi$, so the integral
equals $|x|\cdot2\pi\int_{-1}^1|t|\dd t=2\pi|x|$. In general $c_d>0$ by the
same computation.
\end{proof}

\begin{lemma}[Strict conditional negative definiteness]\label{lem:cnd}
If $\sigma(\R^d)=0$, then $I(\sigma)\le0$, with equality iff $\sigma=0$.
\end{lemma}
\begin{proof}
\emph{Dimension one.} Let $\tau$ be such a measure on $\R$ and put
$T(u)=\tau((-\infty,u])$. For real $s,t$,
\[
|s-t|=\int_\R\big(\mathbf 1\{s\le u<t\}+\mathbf 1\{t\le u<s\}\big)\dd u .
\]
Integrate against $\tau\otimes\tau$; Fubini applies since $|\tau|$ is finite
and the support is bounded. Using $\tau(\{s\le u\})=T(u)$ and
$\tau(\{t>u\})=-T(u)$,
\[
I(\tau)=-2\int_\R T(u)^2\dd u\le0 .
\]
\emph{Projections.} For $\omega\in\Sph^{d-1}$ let
$\pi_\omega(x)=\langle x,\omega\rangle$. By Lemma~\ref{lem:crofton} and
Fubini,
\[
I(\sigma)=\frac1{c_d}\int_{\Sph^{d-1}}
I\big((\pi_\omega)_\#\sigma\big)\dd\sigma(\omega)\le0 .
\]
\emph{Equality.} If $I(\sigma)=0$, then the continuity of
$\omega\mapsto(\pi_\omega)_\#\sigma$ forces every projection to have
$T\equiv0$, so every $(\pi_\omega)_\#\sigma$ vanishes. Hence the Fourier
transform $\hat\sigma(r\omega)=\widehat{(\pi_\omega)_\#\sigma}(r)$ vanishes
identically and $\sigma=0$ (Cram\'er--Wold).
\end{proof}

\begin{remark}
Equivalently,
$I(\sigma)=-C_d\int_{\R^d}|\hat\sigma(\xi)|^2|\xi|^{-d-1}\dd\xi$ with $C_d>0$.
This is the energy distance of Sz\'ekely and Rizzo~\cite{SR13}, and in the
language of Nickolas and Wolf it says that $\R^d$ is \emph{strictly
quasihypermetric}.
\end{remark}

\begin{lemma}[Energy identity]\label{lem:identity}
For measures $\mu,\nu$,
$I(\mu)=I(\nu)+2\int U_\nu\dd(\mu-\nu)+I(\mu-\nu)$.
\end{lemma}
\begin{proof}
Expand $I(\nu+(\mu-\nu))$ by bilinearity and use symmetry of the kernel.
\end{proof}

\subsection{Bj\"orck's theorem}

\begin{proof}[Proof of Theorem~\ref{thm:bjorck}]
\emph{Existence.} $\Prob(A)$ is weak-$*$ compact, and $I$ is weak-$*$
continuous because the kernel is continuous on $A\times A$.

\emph{Uniqueness.} If $\mu_0\ne\mu_1$ both attain $\Mp$, then by bilinearity
and Lemma~\ref{lem:cnd},
$I\big(\tfrac12(\mu_0+\mu_1)\big)=\tfrac12I(\mu_0)+\tfrac12I(\mu_1)-\tfrac14I(\mu_0-\mu_1)>\Mp$,
a contradiction.

\emph{Equilibrium conditions.} For $x\in A$ and
$\mu_\varepsilon=(1-\varepsilon)\mu_0+\varepsilon\delta_x$,
\[
I(\mu_\varepsilon)=(1-\varepsilon)^2\Mp+2\varepsilon(1-\varepsilon)U_{\mu_0}(x),
\]
whose right derivative at $\varepsilon=0$, namely $2(U_{\mu_0}(x)-\Mp)$, must
be $\le0$. Since $\int U_{\mu_0}\dd\mu_0=\Mp$ and $U_{\mu_0}\le\Mp$, equality
holds $\mu_0$-a.e., hence on $\supp\mu_0$ by continuity.

\emph{Support on the boundary.} Let $A=K$ be a convex body. The function
$U=U_{\mu_0}$ is convex. Suppose $U(x)=\Mp$ at an interior point $x$. For any
$y\in K$ write $x=ty+(1-t)z$ with $z\in K$ and $t\in(0,1)$. Then
$\Mp=U(x)\le tU(y)+(1-t)U(z)\le\Mp$, so $U(y)=\Mp$. Thus $U$ would be constant
on $K$, in particular on an open set. But in the sense of distributions
$\Delta U=(d-1)\int|x-y|^{-1}\dd\mu_0(y)>0$, so $U$ is not constant on any
open set. Hence $U<\Mp$ on the interior, and $\supp\mu_0\subset\partial K$.
\end{proof}

\subsection{Points versus measures}

\begin{proof}[Proof of Theorem~\ref{thm:gamma}]
Let $m_n$ be the maximum pair average. For points $\xi_1,\dots,\xi_n$ with
empirical measure $\nu_n=\frac1n\sum\delta_{\xi_j}$,
\[
\frac{2}{n(n-1)}\sum_{j<k}|\xi_j-\xi_k|=\frac{n}{n-1}\cdot\frac1{n^2}\sum_{j,k}|\xi_j-\xi_k|
=\frac{n}{n-1}\,I(\nu_n)\le\frac{n}{n-1}\,\Mp,
\]
which gives the upper bound. For the lower bound, let $\xi_1,\dots,\xi_n$ be
independent with law $\mu_0$. Each pair has
$\mathbb{E}|\xi_j-\xi_k|=I(\mu_0)=\Mp$, so the expected pair average equals
$\Mp$, and some configuration attains at least its expectation. Letting
$n\to\infty$ gives $\gamma_1=\Mp$.
\end{proof}

The same computation shows $\sum_{j,k}|\xi_j-\xi_k|\le\Mp(A)\,n^2$ for every
configuration in every compact $A$. This is the inequality of
Corollary~\ref{cor:equalpot}.

\subsection{Polytopes}

\begin{proof}[Proof of Theorem~\ref{thm:polytope}]
Let $\nu$ be the extremal measure of the finite set $V$. By
Theorem~\ref{thm:bjorck}(2), applied to $V$, $U_\nu\le\Mp(V)$ on $V$. Since
$U_\nu$ is convex, it is at most $\Mp(V)$ on $\conv V=K\supseteq A$. For
$\mu\in\Prob(A)$, Lemma~\ref{lem:identity} gives
\[
I(\mu)=I(\nu)+2\Big(\int U_\nu\dd\mu-\int U_\nu\dd\nu\Big)+I(\mu-\nu)
\le\Mp(V)+2\big(\Mp(V)-\Mp(V)\big)+0=\Mp(V),
\]
using $\int U_\nu\dd\nu=I(\nu)=\Mp(V)$ and Lemma~\ref{lem:cnd}. Equality
forces $\mu=\nu$. Since $V\subseteq A$, also $\Mp(A)\ge\Mp(V)$.
\end{proof}

\begin{proof}[Proof of Corollary~\ref{cor:equalpot}]
The measure $\nu=\sum w_i\delta_{v_i}$ is a probability measure on $V$ with
$U_\nu\equiv\lambda$ on $V$. By Proposition~\ref{prop:signed} below it is the
extremal measure of $V$, and $\Mp(V)=\lambda=I(\nu)$. Theorem~\ref{thm:polytope}
transfers this to $A$. The finite-$n$ inequality was proved above. Equality
holds for the stated configuration because its empirical measure is $\nu$.
\end{proof}

The triangle weights are checked directly: the unnormalized potential at $A$
is $w_Bc+w_Cb=bc\big[(c+a-b)+(a+b-c)\big]=2abc$, and similarly at $B$ and $C$.
The normalizing constant is
$\sum a(b+c-a)=2(ab+bc+ca)-(a^2+b^2+c^2)$. For the regular $m$-gon with
circumradius $1$ the distances from a vertex are $2\sin(\pi k/m)$, and
$\sum_{k=0}^{m-1}2\sin\frac{\pi k}{m}=2\cot\frac{\pi}{2m}$.

\subsection{Finite sets and the kite}\label{app:kite}

\begin{proposition}\label{prop:finite}
Let $V=\{v_1,\dots,v_m\}$, $m\ge2$, be distinct points with distance matrix
$D$. Then $D$ is invertible, $s=D^{-1}\mathbf1$ satisfies
$\mathbf1^{\mathsf T}s>0$, and $M(V)=1/\mathbf1^{\mathsf T}s$, attained only
by $w=s/\mathbf1^{\mathsf T}s$.
\end{proposition}
\begin{proof}
\emph{Invertibility.} Suppose $Dx=0$ with $x\ne0$. If
$\mathbf1^{\mathsf T}x=0$, then $x^{\mathsf T}Dx=0$ contradicts
Lemma~\ref{lem:cnd}. Otherwise scale so that $\mathbf1^{\mathsf T}x=1$. For
each $i$, $y=e_i-x$ has mass zero and
$y^{\mathsf T}Dy=D_{ii}-2e_i^{\mathsf T}Dx+x^{\mathsf T}Dx=0$, so $x=e_i$ for
every $i$, which is impossible.

\emph{Sign.} If $\mathbf1^{\mathsf T}s=0$, then $s^{\mathsf T}Ds=0$ forces
$s=0$, which is impossible. If $\mathbf1^{\mathsf T}s<0$, put
$c=-\mathbf1^{\mathsf T}s/m>0$ and $y=s+c\mathbf1$, which has mass zero. Then
$y^{\mathsf T}Dy=\mathbf1^{\mathsf T}s+2cm+c^2\mathbf1^{\mathsf T}D\mathbf1=cm+c^2\mathbf1^{\mathsf T}D\mathbf1>0$,
contradicting Lemma~\ref{lem:cnd}.

\emph{Value.} $w$ has $Dw=\lambda\mathbf1$ with
$\lambda=1/\mathbf1^{\mathsf T}s$; apply Proposition~\ref{prop:signed}.
\end{proof}

\paragraph{The kite formula.} Write the unnormalized weights as $p,p,q,r$
with every potential equal to $1$. With $s=\sqrt{1+\varepsilon^2}$ and
$t=1+\varepsilon$, the equations are
\[
2p+\sqrt2\,q+sr=1,\qquad 2\sqrt2\,p+tr=1,\qquad 2sp+tq=1 .
\]
Eliminating $q$ and $r$ gives $p=(\sqrt2+s-t)/(4\sqrt2\,s-2t)$, and
$M(V)=1/(2p+q+r)$ simplifies to the formula of Section~\ref{sec:kite} (using
$s^2=1+\varepsilon^2$). At $\varepsilon=0$ it reduces to $(3+\sqrt2)/4$. The
threshold $\varepsilon^*$ is where the potential of the triangle measure
$(w_A,w_B,w_C)$ at $D$ equals $2/(2\sqrt2-1)$, namely
$\sqrt2\,s+(\sqrt2-1)t=2$; this is equivalent to the displayed equation.

\subsection{Signed equilibria and the criterion}\label{app:signed}

\begin{proposition}[Signed equilibrium]\label{prop:signed}
Suppose $\nu\in\Sgn(A)$ has $U_\nu\equiv\lambda$ on $A$. Then every
$\mu\in\Sgn(A)$ satisfies $I(\mu)=\lambda+I(\mu-\nu)\le\lambda$, so
$M(A)=\lambda$ with unique maximizer $\nu$. Moreover $\Mp(A)=M(A)$ iff
$\nu\ge0$, in which case $\nu=\mu_0$.
\end{proposition}
\begin{proof}
By Lemma~\ref{lem:identity},
$I(\mu)=\lambda+2(\lambda-\lambda)+I(\mu-\nu)$, and Lemma~\ref{lem:cnd}
gives $I(\mu-\nu)\le0$, with equality only for $\mu=\nu$. If $\nu\ge0$ it is
admissible for $\Mp$. If not, $\mu_0\ne\nu$ and $\Mp=I(\mu_0)<\lambda$.
\end{proof}

\begin{lemma}[Uniqueness]\label{lem:signeduniq}
There is at most one signed equilibrium on $A$.
\end{lemma}
\begin{proof}
If $\nu_1,\nu_2$ have constant potentials $\lambda_1,\lambda_2$, then
$\sigma=\nu_1-\nu_2$ has mass zero and
$I(\sigma)=\int(U_{\nu_1}-U_{\nu_2})\dd\sigma=(\lambda_1-\lambda_2)\sigma(A)=0$,
so $\sigma=0$.
\end{proof}

\begin{proof}[Proof of Theorem~\ref{thm:criterion}]
If $U_{\mu_0}$ is constant, Proposition~\ref{prop:signed} gives $M=\Mp$.
Conversely, suppose $\Mp=M$. Then $M$ is finite and $\mu_0$ maximizes $I$ over
$\Sgn(A)$. For $x\in A$ and \emph{every} real $\varepsilon$,
$\mu_\varepsilon=\mu_0+\varepsilon(\delta_x-\mu_0)\in\Sgn(A)$ and
\[
I(\mu_\varepsilon)=M+2\varepsilon\big(U_{\mu_0}(x)-M\big)+\varepsilon^2I(\delta_x-\mu_0).
\]
Since $\varepsilon$ may have either sign, maximality forces $U_{\mu_0}(x)=M$.
\end{proof}

\begin{corollary}\label{cor:gap}
If a signed equilibrium $\nu$ on $A$ exists and is not a positive measure,
then $\Mp(A)<M(A)$ and $\supp\mu_0\ne A$.
\end{corollary}
\begin{proof}
The inequality is Proposition~\ref{prop:signed}. If $\supp\mu_0=A$, then
$U_{\mu_0}\equiv\Mp$ on $A$ by Theorem~\ref{thm:bjorck}, so $\mu_0$ is a
signed equilibrium and equals $\nu$ by Lemma~\ref{lem:signeduniq}. This is a
contradiction.
\end{proof}

\begin{proof}[Proof of Proposition~\ref{prop:bracket}]
Apply Lemma~\ref{lem:identity} with $\mu=\mu_0$:
$\Mp=2\int U_\nu\dd\mu_0-I(\nu)+I(\mu_0-\nu)$. Here
$\int U_\nu\dd\mu_0\le\max_AU_\nu$ because $\mu_0\ge0$ has mass one, and
$I(\mu_0-\nu)\le0$.
\end{proof}

\subsection{Hollow polygons}

The following lemma gives a direct proof of Corollary~\ref{cor:wilson}(1),
without appeal to Wilson's theorem.

\begin{lemma}\label{lem:polygon}
Let $A\subset\R^2$ be compact and contain two non-parallel line segments,
for instance the boundary of a convex polygon. Then no probability measure on
$A$ has constant potential on $A$. Consequently $\Mp(A)<M(A)$.
\end{lemma}
\begin{proof}
	Suppose $U_\mu\equiv c$ on $A$, and let $e\subset A$ be a segment, traversed
	linearly as $x(t)$. For each $y$, the function $f_y(t)=|x(t)-y|$ is convex,
	and $f_y''=0$ on the interior of $e$ implies that $y$ lies on the supporting
	line $\ell_e$ of $e$. (If $y\notin\ell_e$, then
	$f_y''=|x'|^2\operatorname{dist}(y,\ell_e)^2/|x-y|^3>0$; if $y$ lies in the
	interior of $e$, then $f_y''$ has an atom at $y$.) The second distributional
	derivatives $f_y''$ are nonnegative measures whose $\mu$-average is
	$(U_\mu\circ x)''=0$ on the interior of $e$. Testing against a countable
	dense family of nonnegative test functions shows that $f_y''=0$ there for
	$\mu$-a.e.\ $y$, so $\supp\mu\subset\ell_e$. Doing the same for a second
	segment $e'$ not parallel to $e$ gives
	$\supp\mu\subset\ell_e\cap\ell_{e'}$, a single point $p$. So $\mu=\delta_p$.
	But $|x-p|$ is not constant along $e$, a contradiction. The inequality
	follows from Theorem~\ref{thm:criterion}.
\end{proof}

\subsection{The width bound}\label{app:width}

\begin{proof}[Proof of Proposition~\ref{prop:width}]
	\emph{Segments.} Let $\tau$ be a signed mass-one measure on an interval
	$[\alpha,\beta]$ of length $L$, and put $T(u)=\tau((-\infty,u])$. The
	computation in the proof of Lemma~\ref{lem:cnd}, now with
	$\tau(\{t>u\})=1-T(u)$, gives
	\[
	I(\tau)=2\int_\alpha^\beta T(u)\big(1-T(u)\big)\dd u\ \le\ 2\cdot\tfrac14\cdot L=\tfrac L2,
	\]
	since $T(1-T)\le\frac14$ for every real $T$. Equality forces $T\equiv\frac12$
	on $(\alpha,\beta)$, i.e.,\ $\tau=\tfrac12(\delta_\alpha+\delta_\beta)$.
	
	\emph{The plane.} For $\nu\in\Sgn(X)$, Lemma~\ref{lem:crofton} (with $c_2=4$)
	and Fubini give
	\[
	I(\nu)=\frac14\int_0^{2\pi}I\big((\pi_\theta)_\#\nu\big)\dd\theta .
	\]
	The projected measure lives on an interval of length $w(\theta)$, the width
	of $\conv X$ in direction $\theta$, so
	\[
	I(\nu)\le\frac18\int_0^{2\pi}w(\theta)\dd\theta=\frac18\cdot2P=\frac P4 .
	\]
	Here $w(\theta)=h(\theta)+h(\theta+\pi)$, where $h$ is the support function,
	and Cauchy's formula gives $\int_0^{2\pi}h\dd\theta=P$.
	
	\emph{Higher dimensions.} The same argument with
	$|x|=b(n)\int_{\Sph^{n-1}}|\langle x,t\rangle|\dd\lambda(t)$, where
	$\lambda$ is the normalized surface measure, gives
	$I(\nu)\le\tfrac12\,b(n)\int w\dd\lambda=\tfrac12\,b(n)\,\bar w$.
	
	\emph{Equality.} For a segment, equality holds by the first part. For the
	solid ball, it holds because $\bar w=2$ and $M(B^n)=b(n)$~\cite{HNW11}.
\end{proof}

\section{Rendezvous numbers}\label{app:rendezvous}

Gross~\cite{Gross64} considered the two-person zero-sum game on a compact
connected metric space $X$. Each player picks a point, and the payoff is
their distance. By Glicksberg's minimax theorem~\cite{Glicksberg50} the game
has a value $r(X)$ and optimal mixed strategies, i.e.,\ probability measures:
\[
r(X)=\max_{\mu\in\Prob(X)}\min_X U_\mu=\min_{\nu\in\Prob(X)}\max_X U_\nu .
\]
Gross called $r(X)$ the \emph{rendezvous value}. He stated his theorem for
equal-weight finite sets, and Elton generalized it to measures
(see~\cite{CMY86}): every $\mu\in\Prob(X)$ has $U_\mu(P)=r(X)$ at some
$P\in X$, and $r(X)$ is the unique constant with this property. Gross also
proved $D/2\le r<D$, where $D$ is the diameter, and asked for the sharp upper
bound of $r/D$ among plane sets.

\begin{theorem}[Wolf~\cite{Wolf97}]\label{thm:rendezvous}
For compact connected $X\subset\R^d$, $r(X)\le\Mp(X)$, with equality iff
$U_{\mu_0}$ is constant on $X$.
\end{theorem}
\begin{proof}
First, $r\le\max U_{\mu_0}=\Mp$ by Theorem~\ref{thm:bjorck}. If
$U_{\mu_0}\equiv\Mp$, then for every $\mu$ we have
$\int U_\mu\dd\mu_0=\int U_{\mu_0}\dd\mu=\Mp$, so
$\min U_\mu\le\Mp\le\max U_\mu$, and uniqueness gives $r=\Mp$. Conversely, if
$r=\Mp$, let $\mu_1$ be optimal for the maximizing player. Then
$U_{\mu_1}\ge\Mp$ on $X$, so $I(\mu_1)\ge\Mp$, whence $\mu_1=\mu_0$ and
$U_{\mu_0}\ge\Mp$. Combined with $U_{\mu_0}\le\Mp$, this gives
$U_{\mu_0}\equiv\Mp$.
\end{proof}

Combined with Theorem~\ref{thm:criterion}, for compact connected
$X\subset\R^d$:
\[
r(X)=\Mp(X)\iff \Mp(X)=M(X).
\]
So positivity also decides whether the rendezvous value equals the maximal
energy. It follows that $r=\Mp=M$ for circles, spheres and segments, and
numerically for $E_{2,1}$. Wolf remarked that $r$ was unknown for a general
ellipse~\cite{Wolf97}; for $E_{2,1}$ it equals the 50-digit value of
Section~\ref{sec:E21}, conditional on the positivity of $\rho$. By contrast,
$r<\Mp$ for $E_{3,1}$, $E_{5,1}$, $\Lell$, every hollow convex polygon and
every convex body. For the unit disk, for instance, $U_{\mu_0}$ at the centre
equals $1<4/\pi$.

\section{Numerical methods and data}\label{app:numerics}

All computations were done in Mathematica; dense linear algebra was
multithreaded through MKL.

\subsection{What is rigorous}\label{app:rigor}

\begin{itemize}
	\item \textbf{Lower bounds for $\Mp$.} For any nonnegative weights $w$ on any
	node set, normalized to total mass one, $I(w)\le\Mp$. This holds whatever
	the accuracy of the solver that produced $w$, so these bounds are rigorous
	once $I(w)$ is evaluated with control of rounding (e.g.\ interval
	arithmetic).
	\item \textbf{Lower bounds for $M$.} Likewise $I(w)\le M$ for \emph{signed}
	weights of total mass one, and the Galerkin values satisfy
	$\lambda_K=I(\nu_K)\le M$, nondecreasing in $K$. We have no rigorous
	\emph{upper} bound yet for $M$; the printed digits of $M$ are those on which
	successive truncations agree.
	\item \textbf{Upper bounds for $\Mp$} come from Proposition~\ref{prop:bracket}
	and require the maximum of $U_\nu$ over the whole curve or surface. In our
	computations it was found by grid search and local maximization, which is
	not rigorous. Lemma~\ref{lem:curv} shows how to make it so.
	\item \textbf{Signs of $\rho$} are numerical; a rigorous version would need a
	truncation-error bound smaller than $|\min\rho|$.
	\item \textbf{Supports.} Dropping nodes with negative unconstrained weight is
	a heuristic. The reported supports are those of the computed discrete
	optima, characterized by the Karush--Kuhn--Tucker conditions.
\end{itemize}

\begin{lemma}[Curvature correction]\label{lem:curv}
	Let $x:[\alpha,\beta]\to\R^d$ be $C^2$ with $|x''|\le\kappa$, and let
	$\mu\in\Prob(\R^d)$ have compact support. Then $g(t)=U_\mu(x(t))$ satisfies
	$g''\ge-\kappa$ in the sense of distributions. Hence
	\[
	\max_{[\alpha,\beta]}g\ \le\ \max\{g(\alpha),g(\beta)\}+\tfrac{\kappa}{8}(\beta-\alpha)^2 .
	\]
\end{lemma}
\begin{proof}
	For $f(t)=|x(t)-v|$ and $x(t)\ne v$,
	\[
	f''=\frac{|x'|^2|x-v|^2-\big((x-v)\cdot x'\big)^2}{|x-v|^3}
	+\frac{(x-v)\cdot x''}{|x-v|}\ \ge\ -|x''|,
	\]
	by Cauchy--Schwarz; where $x(t)=v$, $f$ has a convex kink. Averaging over
	$\mu\ge0$ gives $g''\ge-\kappa$. Thus $g(t)+\frac\kappa2(t-c)^2$, with $c$
	the midpoint, is convex, and its maximum on $[\alpha,\beta]$ is attained at
	an endpoint.
\end{proof}

For the ellipse $\kappa=a$. With the grids used here, of spacing
$2\pi/(40n)$, the correction $\kappa h^2/8$ is below $2\times10^{-11}$ at
$n=30000$, comparable to the floating-point error in the $n$-term sums. After
symmetrizing the weights over the four symmetries of the ellipse, which
keeps a probability measure and makes the first-quadrant search exact, the
ellipse brackets thus become rigorous with upper ends raised by less than
$10^{-10}$. Combined with the energy of the unconstrained signed weights on
the same nodes, which is a lower bound for $M$, this would prove
$\Mp<M$ for $E_{3,1}$ and $E_{5,1}$. We have not carried out this
verification. For a surface parametrized over a rectangle of parameters, the
same argument gives a correction $\kappa\delta^2/8$ on grid cells of
diameter $\delta$, where $|D^2x[h,h]|\le\kappa|h|^2$. For the ellipsoid in
coordinates $(\theta,\varphi)$ one may take $\kappa=2\max(a,b,c)$; with the
$200\times200$ octant grid used here, the correction (about $4\times10^{-5}$)
is too large, and a grid roughly ten times finer would be needed.

\paragraph{A test case.} For $n$ equally weighted nodes at the $n$th roots of
unity on the unit circle, set $x=\pi/(2n)$. The potential is
$\frac2n\cot x$ at the nodes and $\frac2n\csc x$ midway between them, so the
bracket has width $\frac4n\tan\frac\pi{4n}\approx\pi/n^2$. The lower and upper
errors are $\pi/(3n^2)$ and $2\pi/(3n^2)$. For $n=800$ the width is
$4.90874\times10^{-6}$, which the code reproduces digit for digit.

\subsection{Ellipses: dense discretization}\label{app:dense}

Take $n$ nodes uniform in $t$, solve $Ds=\mathbf1$, remove nodes with
$s_i<0$, and repeat until all weights are positive. The resulting
$w=s/\sum s_i$ gives the lower bound $I(w)=1/\sum s_i$, and
Proposition~\ref{prop:bracket} gives the upper bound. Both errors decay like
$n^{-2}$, so Richardson extrapolation $R=\tfrac13(4X(2n)-X(n))$ is very
effective.

\paragraph{$E_{2,1}$.} We first report values for $E_{1,1/2}$, which equal
$\tfrac12$ of those for $E_{2,1}$. Every node was kept at every $n$ tested.
The errors are approximately $1.221/n^2$ (lower bound) and $3.913/n^2$ (upper
bound). Richardson values:
\begin{center}
\begin{tabular}{@{}lcc@{}}
\toprule
$(n,2n)$ & from lower bounds & from upper bounds\\
\midrule
400/800 & 1.054027106093 & 1.05402710765\\
1600/3200 & 1.054027106093 & 1.05402710609\\
6400/12800 & 1.054027106087 & 1.05402710609\\
12800/25600 & 1.054027106093 & 1.05402710609\\
\bottomrule
\end{tabular}
\end{center}
The predicted discrete optimum on LSW's grid, $2(\Mp(E_{1,1/2})-1.2202/N^2)$
with $N=4095$, is $2.108054066650$; LSW reported $2.108054066642$.
For $E_{2,1}$ itself, runs up to $n=35000$ keep every node; at $n=35000$,
$2.10805421018\le\Mp\le2.10805421856$, and Richardson extrapolation of the
$16000/32000$ pair gives $2.10805421217$ (lower) and $2.10805421218$ (upper).

\paragraph{$E_{3,1}$.}
\begin{center}
	\begin{tabular}{@{}rllll@{}}
		\toprule
		$n$ & lower & upper & width & support\\
		\midrule
		400 & 3.06770258391 & 3.06780511818 & $1.03\times10^{-4}$ & 286\\
		1600 & 3.06772567409 & 3.06773208532 & $6.41\times10^{-6}$ & 1146\\
		6400 & 3.0677271172 & 3.06772751791 & $4.01\times10^{-7}$ & 4578\\
		12800 & 3.06772718935 & 3.06772728953 & $1.00\times10^{-7}$ & 9158\\
		20000 & 3.06772720356 & 3.06772724459 & $4.10\times10^{-8}$ & 14306\\
		30000 & 3.06772720903 & 3.06772722727 & $1.82\times10^{-8}$ & 21458\\
		\bottomrule
	\end{tabular}
\end{center}
The width satisfies $\text{width}\cdot n^2\approx16.4$ for $n\ge12800$.
Richardson extrapolation of the $15000/30000$ pair gives $3.06772721341$ from
both the lower and the upper bounds, and the $14000/28000$ pair gives the
same value. The upper bound at $n=30000$ lies $4.4\times10^{-8}$ below $M$.
About $71.5\%$ of the nodes are kept at every $n$.

\paragraph{$E_{5,1}$.}
\begin{center}
\begin{tabular}{@{}rllll@{}}
\toprule
$n$ & lower & upper & width & support\\
\midrule
400 & 5.03937350639 & 5.03954932517 & $1.76\times10^{-4}$ & 162\\
1600 & 5.03941330497 & 5.03942430774 & $1.10\times10^{-5}$ & 646\\
6400 & 5.03941579224 & 5.03941647997 & $6.88\times10^{-7}$ & 2578\\
12800 & 5.0394159166 & 5.03941608854 & $1.72\times10^{-7}$ & 5158\\
20000 & 5.03941594108 & 5.0394160115 & $7.04\times10^{-8}$ & 8058\\
30000 & 5.03941595051 & 5.03941598181 & $3.13\times10^{-8}$ & 12082\\
\bottomrule
\end{tabular}
\end{center}
The width satisfies $\text{width}\cdot n^2\approx28.2$ for $n\ge12800$.
Richardson extrapolation of the $15000/30000$ pair gives $5.039415958057$ from
both the lower and the upper bounds. About $40.3\%$ of the nodes are kept at
every $n$.

\subsection{Ellipses: spectral Galerkin}\label{app:spectral-ellipse}

\paragraph{Kernel.} Let $A_0=\frac{a^2+b^2}2$, $B_0=\frac{a^2-b^2}2$ and
$q=\frac{a-b}{a+b}$. Then
\[
|x(t)-x(s)|=\left|\sin\tfrac{t-s}2\right|g(t+s),\qquad
g(u)=2\sqrt{A_0-B_0\cos u}=(a+b)\,|1-qe^{iu}| .
\]
The first factor has Fourier coefficients $\alpha_m=-2/(\pi(4m^2-1))$. The
second has coefficients $\beta_p=\beta_{-p}$, which are given exactly by
\[
\beta_p=(a+b)\,\frac{(-\tfrac12)_p}{p!}\,q^p\;{}_2F_1\!\left(-\tfrac12,\,p-\tfrac12;\,p+1;\,q^2\right),\qquad p\ge0,
\]
obtained by expanding $(1-qe^{\pm iu})^{1/2}$. As a check,
$\beta_0=(a+b)\,{}_2F_1(-\frac12,-\frac12;1;q^2)=P/\pi$ is the Gauss--Kummer
series for the perimeter. For $p\ge1$ every term of the ${}_2F_1$ after the
first is negative and the series sum is positive, so $\beta_p$ is computed
with full relative precision. The decay is $\beta_p\sim q^p$.

\paragraph{System.} By symmetry under $t\mapsto-t$ and $t\mapsto t+\pi$,
$\rho(s)=\sum_J\rho_{2J}e^{2iJs}$ with $\rho_{2J}=\rho_{-2J}$, and the
potential has coefficients
$U_{2Q}=2\pi\sum_J\alpha_{Q+J}\beta_{Q-J}\rho_{2J}$. Impose $U_{2Q}=0$ for
$0<|Q|\le K$ and $\rho_0=1/(2\pi)$; then $\lambda_K=U_0$.

\paragraph{A precision lesson.} Computing $\beta_p$ by numerical quadrature
sums $O(1)$ terms that cancel down to about $q^p$. In significance
arithmetic this destroys all digits by $K=80$. The closed form avoids the
problem entirely.

\paragraph{Results.} For $E_{2,1}$ ($q=\frac13$) the error in $\lambda_K$ is
about $1.35\times10^{-3}\cdot3^{-2K}$ for $E_{1,1/2}$. For $E_{5,1}$
($q=\frac23$), working with $\max(60,K+30)$ digits:
\begin{center}
\small
\begin{tabular}{@{}rlll@{}}
\toprule
$K$ & $\lambda_K$ & $|\rho_{2K}|$ & $\min_t\rho$\\
\midrule
10 & 5.0394115880284988929057069419238219419569932642585 & $6.36\times10^{-3}$ & $-0.006466425799$\\
20 & 5.0394162861805010595144883985910747933481747916637 & $1.23\times10^{-4}$ & $-0.0005165414917$\\
40 & 5.0394162866751156051838201513081845136263829227462 & $4.15\times10^{-8}$ & $-0.0004111309058$\\
80 & 5.0394162866751156203052814396774491962128628511346 & $4.24\times10^{-15}$ & $-0.0004110980116$\\
160 & 5.0394162866751156203052814396774901490677251309166 & $3.94\times10^{-29}$ & $-0.0004110980116$\\
\bottomrule
\end{tabular}
\end{center}
Successive differences of $\lambda_K$ are $4.9\times10^{-10}$,
$1.5\times10^{-17}$ and $4.1\times10^{-32}$.
For $E_{3,1}$ ($q=\tfrac12$) the same code gives $\lambda_{40}$, $\lambda_{80}$
and $\lambda_{160}$ with tails $5.3\times10^{-13}$, $5.4\times10^{-25}$ and
$5.1\times10^{-49}$; $\lambda_{80}$ and $\lambda_{160}$ agree to 50 digits.
The minimum of $\rho$ over a grid of spacing $\pi/800$ is
$-0.0005105319482$ for $K\ge80$.

\subsection{The ellipsoid: spectral Galerkin}\label{app:spectral-ellipsoid}

\paragraph{Basis.} Uniqueness and the reflection symmetries make $\rho$ even
in each coordinate. We use the real orthonormal spherical harmonics with $l$
even, $m$ even, of cosine type; there are about $L^2/8$ of them up to degree
$L$.

\paragraph{Crofton plus Funk--Hecke.} Let $A=\operatorname{diag}(a,b,c)$ and
$n=A\omega/|A\omega|$. Lemma~\ref{lem:crofton} gives
\[
|A(u-v)|=\frac1{2\pi}\int_{\Sph^2}|A\omega|\,
\big|\langle u,n\rangle-\langle v,n\rangle\big|\dd\sigma(\omega),
\]
so the energy is an average of one-dimensional energies of projections. By
Funk--Hecke, the projection of $Y_{lm}(u)\dd\sigma(u)$ onto the axis $n$ has
density $2\pi Y_{lm}(n)P_l(t)$. The singular kernel thus enters only through
the exact rationals
\[
e_{ll}=\frac{-8}{(2l-1)(2l+1)(2l+3)},\qquad
e_{l,l+2}=\frac{4}{(2l+1)(2l+3)(2l+5)},\qquad e_{ll'}=0\ (|l-l'|>2),
\]
where $e_{ll'}=\iint|s-t|P_l(s)P_{l'}(t)\dd s\dd t$. Moreover,
$F_l(s)=\int_{-1}^1|s-t|P_l(t)\dd t$ equals $s^2+1$ for $l=0$ and
$6(1-s^2)^2C^{(5/2)}_{l-2}(s)/((l-1)l(l+1)(l+2))$ for even $l\ge2$. This form
avoids cancellation when the potential is evaluated.

\paragraph{System.} For $\rho=\sum_ic_iY_i$ we have $I(\rho)=c^{\mathsf T}Gc$
with
\[
G_{ij}=2\pi\,e_{l_il_j}\int_{\Sph^2}|A\omega|\,Y_i(n)Y_j(n)\dd\sigma(\omega),
\]
computed by a product Gauss rule on one octant. Solve $Gx=e_{00}$; then
$\lambda_L=1/(4\pi x_{00})$ and $c=x/(\sqrt{4\pi}\,x_{00})$. For $A=I$, $G$ is
diagonal with the sphere eigenvalues $-16\pi/((2l-1)(2l+1)(2l+3))$, and
$\lambda_0$ is the energy of the uniform measure. The ratio of the $l=0$ and
$l=2$ eigenvalues is $35$; this is the amplification mentioned in
Section~\ref{sec:ellipsoid}.

\paragraph{Checks (machine precision).} The basis has Gram error
$1.5\times10^{-13}$. Spheres of radius $1$ and $2$ give $4/3$ and $8/3$ to
$4\times10^{-16}$. Doubling the extra quadrature nodes changes $\lambda_{24}$
by $9.5\times10^{-15}$, and the potential is flat to $2\times10^{-15}$. The
convergence is:
\begin{center}
\begin{tabular}{@{}rlll@{}}
\toprule
$L$ & $\lambda_L$ & $|c_L|$ & $\min\rho$\\
\midrule
0 & 1.488709959082897 & $2.82\times10^{-1}$ & $+0.0795775$\\
4 & 1.511268670191912 & $4.07\times10^{-2}$ & $-0.0127418$\\
8 & 1.511269672970254 & $8.56\times10^{-4}$ & $-0.0168489$\\
16 & 1.511269673044751 & $1.96\times10^{-7}$ & $-0.0169305$\\
24 & 1.511269673044731 & $3.29\times10^{-11}$ & $-0.0169305$\\
\bottomrule
\end{tabular}
\end{center}
The error in $\lambda_L$ is about $|c_L|^2$, as expected for a variational
method.

\paragraph{High precision.} With default significance arithmetic,
Mathematica reports fewer digits as $L$ grows. The cause is that the smallest
entries of $G$ are about $10^{-0.65L}$ times the largest, and the
\emph{relative} precision of the whole matrix is set by them (see the table
below). Their absolute accuracy remains excellent. Resetting the precision of
$G$ before the solve, and checking against a run at higher working precision,
removes the problem.
\begin{center}
\begin{tabular}{@{}rcccc@{}}
\toprule
$L$ & Precision$[G]$ & $\log_{10}(\max|G|/\min|G|)$ & $\log_{10}\kappa(G)$ & Precision$[x_1]$\\
\midrule
16 & 38.11 & 13.08 & 4.16 & 38.11\\
24 & 32.25 & 18.38 & 4.68 & 32.25\\
32 & 26.27 & 23.57 & 5.05 & 26.27\\
40 & 20.25 & 28.69 & 5.34 & 20.25\\
\bottomrule
\end{tabular}
\end{center}
At 60-digit working precision:
\begin{center}
\footnotesize
\begin{tabular}{@{}rl@{}}
\toprule
$L$ & $\lambda_L$\\
\midrule
16 & 1.511269673044740064321404561395408945549724532855818780\\
24 & 1.511269673044740064891740389750224876779268467493957467\\
32 & 1.511269673044740064891740394794742980288877390170950686\\
40 & 1.511269673044740064891740394794743027630511298004227564\\
48 & 1.511269673044740064891740394794743027630511756980212746\\
56 & 1.511269673044740064891740394794743027630511756980217288\\
\bottomrule
\end{tabular}
\end{center}
Runs at 70 digits with more quadrature nodes agree with the 50-digit runs to
within $10^{-50}$ for $L=40$ and $48$.

\paragraph{Point checks.} Solving $Ds=\mathbf1$ without sign constraints on
$1000$ and $2000$ Fibonacci points of the unit sphere gives all weights
positive and energies just below $4/3$. On $\Lell$ the same solves give about
$10\%$ negative weights, consistent with the sign change of $\rho$.

\subsection{The ellipsoid: discrete optimization}\label{app:qp}

\paragraph{Symmetry.} The extremal measure is invariant under the eight sign
changes. Place $n$ Fibonacci nodes $x_i$ in the open first octant, and give
each node and its seven mirror images weight $w_i/8$. The energy becomes
$w^{\mathsf T}\tilde Dw$ with
$\tilde D_{ij}=\frac18\sum_{\epsilon\in\{\pm1\}^3}|x_i-\epsilon x_j|$.

\paragraph{Convexity.} Lemma~\ref{lem:cnd} makes $w\mapsto w^{\mathsf T}\tilde Dw$
strictly concave on the simplex, so its maximizer is unique and characterized 
by the KKT conditions ($U_i=M_n$ on the support, $U_i\le M_n$ elsewhere). On the
simplex $w^{\mathsf T}Jw=1$, where $J$ is the all-ones matrix, so the problem
is the convex quadratic program $\min w^{\mathsf T}(tJ-\tilde D)w$. The
matrix $tJ-\tilde D$ is positive semidefinite once $t$ is at least the signed
discrete constant $1/(\mathbf1^{\mathsf T}\tilde D^{-1}\mathbf1)$. The QP
solution is polished by an exact active-set solve, and the KKT violation is
checked; it was $\lesssim2\times10^{-14}$ in every run.

\paragraph{Results.}
\begin{center}
\begin{tabular}{@{}rrll@{}}
\toprule
octant nodes & support & lower & upper\\
\midrule
501 & 385 & 1.5112542239486 & 1.5113035953411\\
2002 & 1516 & 1.5112593079542 & 1.511266657202\\
8000 & 6027 & 1.5112599272211 & 1.5112608107351\\
16001 & 12029 & 1.5112599834297 & 1.5112602899604\\
20003 & 15038 & 1.5112599921891 & 1.5112602054275\\
25002 & 18781 & 1.5112599984628 & 1.5112601522271\\
30001 & 22534 & 1.5112600022317 & 1.5112601264508\\
\bottomrule
\end{tabular}
\end{center}
Both bounds converge like $n^{-3/2}$. Extrapolations from the triples
$(20003,25002,30001)$ and $(8000,16001,30001)$ agree:
$\Mp(\Lell)\approx1.511260014$. The support fraction tends to about $0.749$.
Two circular caps covering $25.1\%$ of the parameter sphere have angular
radius $\arccos0.749\approx41.5^\circ$. Even $n=501$, i.e.,\ $4008$ effective
points, beats LSW's $1.5112509356$. These upper bounds are not rigorous; 
see Appendix~\ref{app:rigor}.

\section{Mathematica code for ellipses}\label{app:code}

Two programs treat the curve $E_{a,b}$; they solve different problems and
are kept separate.
\begin{itemize}
	\item \texttt{ellipseMplusBracket[a,b,n]} places $n$ point masses uniformly
	in $t$, computes the optimal \emph{nonnegative} weights by an active-set
	solve (Appendix~\ref{app:dense}), and returns the bracket
	$I(w)\le\Mp\le2\max U_w-I(w)$ of Proposition~\ref{prop:bracket}. Its error
	is $O(n^{-2})$. The integer $n$ must be divisible by $4$, so that the
	node set has the symmetries of the ellipse and the maximum of $U_w$ may be
	sought in one quadrant.
	\item \texttt{ellipseMSpectral[a,b,K]} solves the \emph{signed} problem by
	the Fourier--Galerkin method of Appendix~\ref{app:spectral-ellipse}, with
	modes $\cos2Jt$, $0\le J\le K$. It returns $\lambda_K$ (key \texttt{"M"}),
	the coefficients of $\rho$, and $\min\rho$ on a grid. Its error decays like
	$q^{2K}$, $q=\frac{a-b}{a+b}$.
\end{itemize}
The parameters $n$ and $K$ are unrelated: $n$ counts point masses and $K$
counts Fourier modes. The two programs agree when $\rho>0$, since then
$\Mp=M$. When $\rho$ changes sign, the upper end of the bracket for $\Mp$
eventually falls below $M$. As explained in Appendix~\ref{app:rigor},
\texttt{"upper"} is computed by grid search and local maximization and is
not rigorous as it stands.

\subsection{The bracket for $\Mp$}

\begin{lstlisting}
ClearAll[ellipseMplusBracket];
ellipseMplusBracket[a_?NumericQ, b_?NumericQ,
n_Integer?(Mod[#, 4] == 0 &), refine_Integer : 40] :=
Module[{ts, P, D0, active, sol, w, Pact, PT, lower, U, Ut,
tg, G, ug, cand, h, umax},
(* n nodes, uniform in the parameter t *)
ts = N[2 Pi Range[0, n - 1]/n];
P = Transpose[{a Cos[ts], b Sin[ts]}];
D0 = DistanceMatrix[P];
(* active-set solve of D w = const with w > 0 *)
active = Range[n];
While[True,
sol = LinearSolve[D0[[active, active]],
ConstantArray[1., Length[active]]];
If[Min[sol] > 0, Break[]];
active = Pick[active, Positive[sol]]];
w = sol/Total[sol];
Pact = P[[active]];
lower = w . D0[[active, active]] . w;    (* I(w) <= M+ *)
(* potential of the discrete measure *)
PT = Transpose[Pact];
U[{x_?NumericQ, y_?NumericQ}] :=
w . Sqrt[(PT[[1]] - x)^2 + (PT[[2]] - y)^2];
Ut[t_?NumericQ] := U[{a Cos[t], b Sin[t]}];
(* fine grid on the first quadrant: refine points per node gap *)
tg = N[(Pi/2) Range[0, refine n/4]/(refine n/4)];
G = Transpose[{a Cos[tg], b Sin[tg]}];
ug = Join @@ (DistanceMatrix[#, Pact] . w & /@
Partition[G, UpTo[2000]]);
(* polish the 20 largest grid values by local maximization *)
h = 2 Pi/(refine n);
cand = Ordering[ug, -Min[20, Length[ug]]];
umax = Max[ug, Quiet@Table[First@FindMaximum[Ut[t],
{t, tg[[i]] - h, tg[[i]] + h}], {i, cand}]];
<|"n" -> n, "lower" -> lower, "upper" -> 2 umax - lower,
"gap" -> 2 umax - 2 lower, "support" -> Length[active],
"minWeight" -> Min[w], "nodes" -> ts[[active]],
"weights" -> w|>]
\end{lstlisting}

\subsection{The signed constant $M$}

\begin{lstlisting}
ClearAll[ellipseMSpectral];
ellipseMSpectral[a0_, b0_, K_Integer, prec_ : Automatic] :=
Module[{a, b, p0, q, q2, betaL, beta, alpha, Js, Mat, sys, rhs,
x, c, lam, tt},
{a, b} = Rationalize[{a0, b0}, 0];          (* exact inputs *)
p0 = If[prec === Automatic, Max[60, K + 30], prec];
q = (a - b)/(a + b); q2 = N[q^2, p0];
(* exact Fourier coefficients of g(u) = (a+b)|1 - q e^{iu}| *)
betaL = Table[N[(a + b) Pochhammer[-1/2, p]/p! *
If[p == 0, 1, q^p], p0] *
Hypergeometric2F1[-1/2, p - 1/2, p + 1, q2], {p, 0, 2 K}];
beta[p_] := betaL[[Abs[p] + 1]];
alpha[m_] := N[-2/(Pi (4 m^2 - 1)), p0];
Js = Range[-K, K];
Mat = Table[2 Pi alpha[J + Q] beta[Q - J], {Q, Js}, {J, Js}];
(* U_{2Q} = 0 for Q != 0, and normalization rho_0 = 1/(2 Pi) *)
sys = ReplacePart[Mat, K + 1 -> UnitVector[2 K + 1, K + 1]];
rhs = N[UnitVector[2 K + 1, K + 1]/(2 Pi), p0];
x = LinearSolve[sys, rhs];
lam = Mat[[K + 1]] . x;                       (* lambda_K = U_0 *)
c = x[[K + 1 ;;]];                            (* rho_0, ..., rho_2K *)
tt = N[Pi Range[0, 400]/800, p0];             (* grid on [0, Pi/2] *)
<|"M" -> lam, "coeffs" -> c, "tail" -> Abs[Last[c]],
"minRho" -> Min[c[[1]] +
2 Sum[c[[J + 1]] Cos[2 J tt], {J, 1, K}]],
"precision" -> Precision[lam]|>]
\end{lstlisting}

\subsection{Comparing the two}

\begin{lstlisting}
ClearAll[ellipseCompare];
ellipseCompare[a_, b_, n_Integer, K_Integer : 80] :=
Module[{d = ellipseMplusBracket[a, b, n],
s = ellipseMSpectral[a, b, K]},
<|"M+ lower" -> NumberForm[d["lower"], 12],
"M+ upper" -> NumberForm[d["upper"], 12],
"M" -> N[s["M"], 20], "minRho" -> N[s["minRho"], 10],
"support fraction" -> N[d["support"]/n],
"upper < M" -> d["upper"] < N[s["M"]]|>]

(* checks: the circle gives 4/Pi in both programs *)
ellipseMSpectral[1, 1, 8]["M"] - 4/Pi
ellipseMplusBracket[1, 1, 800]["lower"] - (2/800.) Cot[Pi/1600]

(* E_{2,1}: full support, minRho > 0, upper >= M *)
ellipseCompare[2, 1, 6400]
(* E_{5,1}: gaps, minRho < 0, upper < M once n >= 12800 *)
ellipseCompare[5, 1, 12800]
\end{lstlisting}

\section{Mathematica code for the ellipsoid}\label{app:code-ellipsoid}

The surface is given by its semi-axes, e.g.\ \texttt{ax=\{Sqrt[3/2],1,Sqrt[5]/2\}}
for $\Lell$. As for ellipses, two separate programs solve the two problems.
\begin{itemize}
	\item \texttt{ellipsoidMSpectral[ax,L]} solves the \emph{signed} problem by
	the Galerkin method of Appendix~\ref{app:spectral-ellipsoid}, using the real
	spherical harmonics of even degree $l\le L$, even order and cosine type
	(helper \texttt{evenSH}), a product Gauss rule on one octant
	(\texttt{octantRule}), and the exact Legendre energies $e_{ll'}$
	(\texttt{eLeg}). It returns $\lambda_L$ (key \texttt{"M"}), the
	coefficients of $\rho$, the norms of the degree-$l$ parts, $\min\rho$ on a
	grid, and the Gram error of the quadrature. For many digits use
	\texttt{WorkingPrecision->60} and compare with a run at higher precision
	and with more quadrature nodes.
	\item \texttt{ellipsoidMplusQP[ax,n]} places about $n$ Fibonacci nodes in the
	open first octant, symmetrizes the distance matrix over the eight sign
	changes (\texttt{symDist}), and computes the optimal \emph{nonnegative}
	weights by convex quadratic programming followed by an exact active-set
	polish (Appendix~\ref{app:qp}). Its \texttt{"lower"} is a lower bound for
	$\Mp$. \texttt{ellipsoidMplusUpper} then evaluates
	$2\max U_w-I(w)$ by a grid search on the octant and one local
	maximization; as explained in Appendix~\ref{app:rigor}, this upper bound is
	not rigorous as it stands.
\end{itemize}
The degree $L$ and the node count $n$ are unrelated. Because $\rho$ changes
sign on $\Lell$, the two programs do not agree; instead the computed upper
bound for $\Mp$ lies below $M$.

\subsection{The signed constant $M$}

\begin{lstlisting}
ClearAll[evenSH, octantRule, eLeg, ellipsoidMSpectral];

(* orthonormal real spherical harmonics even in x, y, z
(l even, m even, cosine type) at unit vectors pts;
returns {labels {l,m}, matrix points x basis} *)
evenSH[L_Integer, pts_] :=
Module[{x, y, z, w, cols = <||>, cm, Cm, p0, p1, p2, keys},
{x, y, z} = Transpose[pts];
w = x + I y;                       (* Re[w^m] = sin^m(th) cos(m ph) *)
Do[cm = Sqrt[(2 m + 1)/(4 Pi (2 m)!)] (2 m - 1)!!;
Cm = If[m == 0, 1, Sqrt[2] Re[w^m]];
p2 = cm + 0 z; cols[{m, m}] = p2 Cm;          (* l = m *)
p1 = Sqrt[2 m + 3] z p2;                      (* l = m+1 *)
Do[p0 = Sqrt[(4 l^2 - 1)/(l^2 - m^2)] (z p1 -
Sqrt[((l - 1)^2 - m^2)/(4 (l - 1)^2 - 1)] p2);
If[EvenQ[l], cols[{l, m}] = p0 Cm];
p2 = p1; p1 = p0, {l, m + 2, L}], {m, 0, L, 2}];
keys = Sort[Keys[cols]];
{keys, Transpose[cols[#] & /@ keys]}]

(* product rule on the first octant (weights include the factor 8) *)
octantRule[axes_, L_, ex_, prec_] :=
Module[{nt = L + ex, np = L/2 + ex, gz, gw, ph, om, wts, An, nrm},
{gz, gw} = Most[NIntegrate`GaussRuleData[nt, prec]];
ph = N[(Range[np] - 1/2) Pi/(2 np), prec];
om = Flatten[Outer[{Sqrt[1 - #1^2] Cos[#2],
Sqrt[1 - #1^2] Sin[#2], #1} &, gz, ph], 1];
wts = 8 Flatten[Outer[Times, gw,
ConstantArray[N[Pi/(2 np), prec], np]]];
An = Transpose[N[axes, prec] Transpose[om]];
nrm = Sqrt[Total[An^2, {2}]];
{om, wts, nrm, An/nrm}]

(* e_{l l'} = Integrate[|s-t| P_l(s) P_l'(t), {s,-1,1}, {t,-1,1}] *)
eLeg[l_, lp_] := Which[
l == lp, -8/((2 l - 1) (2 l + 1) (2 l + 3)),
Abs[l - lp] == 2,
With[{k = Min[l, lp]}, 4/((2 k + 1) (2 k + 3) (2 k + 5))],
True, 0];

Options[ellipsoidMSpectral] =
{WorkingPrecision -> MachinePrecision, "ExtraNodes" -> 40};
ellipsoidMSpectral[axes_List, L_Integer?(EvenQ[#] && # >= 0 &),
OptionsPattern[]] :=
Module[{prec = OptionValue[WorkingPrecision],
ex = OptionValue["ExtraNodes"], om, wts, nrm, nn, keys, B, B0,
ls, G, x, lam, coef, gram, ug, rhoU, rhoS, dn},
{om, wts, nrm, nn} = octantRule[axes, L, ex, prec];
{keys, B} = evenSH[L, nn];
ls = keys[[All, 1]];
G = 2 Pi (Transpose[B] . ((wts nrm) B)) N[Outer[eLeg, ls, ls], prec];
(* reset precision: significance arithmetic is overly pessimistic here *)
x = LinearSolve[SetPrecision[G, prec],
N[UnitVector[Length[keys], 1], prec]];
lam = 1/(4 Pi x[[1]]);                         (* lambda_L <= M *)
coef = x/(Sqrt[4 Pi] x[[1]]);                  (* coefficients of rho *)
(* diagnostics, machine precision *)
B0 = Last[evenSH[L, N[om]]];
gram = Max[Abs[Transpose[B0] . (N[wts] B0) -
IdentityMatrix[Length[keys]]]];
ug = N@Flatten[Outer[{Sin[#1] Cos[#2], Sin[#1] Sin[#2], Cos[#1]} &,
Subdivide[0, Pi/2, 60], Subdivide[0, Pi/2, 60]], 1];
rhoU = Last[evenSH[L, ug]] . N[coef];          (* density w.r.t. d sigma(u) *)
rhoS = rhoU/((Times @@ N[axes]) Sqrt[Total[
Transpose[Transpose[ug]/N[axes]]^2, {2}]]);  (* w.r.t. area *)
dn = Table[Norm[Pick[N[coef], ls, l]], {l, 0, L, 2}];
<|"axes" -> axes, "L" -> L, "extra" -> ex, "M" -> lam,
"keys" -> keys, "coeffs" -> coef, "degreeNorms" -> dn,
"minRho" -> Min[rhoU], "maxRho" -> Max[rhoU],
"surfaceDensityMinMax" -> MinMax[rhoS], "gramError" -> gram|>]
\end{lstlisting}

\subsection{Bounds for $\Mp$}

\begin{lstlisting}
ClearAll[octantNodes, octantSigns, symDist, ellipsoidMplusQP,
ellipsoidMplusUpper];

(* Fibonacci lattice with 8n points; keep those in the open first
octant (about n of them) *)
octantNodes[n_Integer] := Module[{m = 8 n, k, z, ph, u},
k = N[Range[0, m - 1] + 1/2]; z = 1 - 2 k/m;
ph = N[Pi (3 - Sqrt[5])] k;
u = Transpose[{Sqrt[1 - z^2] Cos[ph], Sqrt[1 - z^2] Sin[ph], z}];
Select[u, Min[#] > 0 &]]

octantSigns = Tuples[{1., -1.}, 3];
(* symmetrized distances (1/8) Sum_s |p_i - s x_j| *)
symDist[P_, X_] := Mean[Table[
DistanceMatrix[P, Transpose[s Transpose[X]]], {s, octantSigns}]]

ellipsoidMplusQP[axes_, n_Integer, tol_ : 10^-9] :=
Module[{u, X, D, m, one, lamS, t, w, S, s, Mn, U, viol, it = 0},
u = octantNodes[n]; X = Transpose[N[axes] Transpose[u]];
D = symDist[X, X]; m = Length[X]; one = ConstantArray[1., m];
lamS = 1/Total[LinearSolve[D, one]];   (* signed discrete constant *)
t = 1.01 lamS;                         (* makes t J - D positive semidefinite *)
w = QuadraticOptimization[
{2 (t ConstantArray[1., {m, m}] - D), ConstantArray[0., m]},
{IdentityMatrix[m, SparseArray], ConstantArray[0., m]},
{{one}, {-1.}}];
(* polish: exact active-set solve starting from the QP support *)
S = Flatten@Position[w, _?(# > tol Max[w] &)];
While[it++ < 100,
s = LinearSolve[D[[S, S]], ConstantArray[1., Length[S]]];
If[Min[s] < 0, S = Pick[S, Positive[s]]; Continue[]];
Mn = 1/Total[s]; w = ConstantArray[0., m]; w[[S]] = Mn s;
U = D . w;
viol = Complement[Flatten@Position[U - Mn, _?(# > 10^-13 &)], S];
If[viol === {}, Break[]]; S = Union[S, viol]];
<|"n" -> m, "effectivePoints" -> 8 m, "signedDiscrete" -> lamS,
"lower" -> w . D . w, "support" -> Length[S],
"kktViolation" -> Max[U] - Mn, "nodes" -> X, "u" -> u,
"weights" -> w|>]

(* upper bound 2 max U - I(w); U has octant symmetry *)
ellipsoidMplusUpper[res_Association, axes_, ng_ : 200] :=
Module[{X = res["nodes"], w = res["weights"], lo = res["lower"],
th, ph, G, Ug, best, Uf, fm},
{th, ph} = Transpose@Flatten[Outer[List, Subdivide[0., Pi/2, ng],
Subdivide[0., Pi/2, ng]], 1];
G = Transpose[N[axes] {Sin[th] Cos[ph], Sin[th] Sin[ph], Cos[th]}];
Ug = symDist[G, X] . w;
best = Transpose[{th, ph}][[First@Ordering[Ug, -1]]];
Uf[a_?NumericQ, b_?NumericQ] := First[symDist[
{N[axes] {Sin[a] Cos[b], Sin[a] Sin[b], Cos[a]}}, X] . w];
fm = First@FindMaximum[{Uf[a, b], 0 <= a <= Pi/2, 0 <= b <= Pi/2},
{{a, best[[1]]}, {b, best[[2]]}}];
<|"lower" -> lo, "upper" -> 2 Max[fm, Max[Ug]] - lo,
"maxU" -> Max[fm, Max[Ug]]|>]
\end{lstlisting}

\subsection{Checks and examples}

\begin{lstlisting}
ax = {Sqrt[3/2], 1, Sqrt[5]/2};

(* basis and quadrature: about 1e-13 *)
ellipsoidMSpectral[ax, 20]["gramError"]
(* sphere of radius 1: about 1e-16 *)
ellipsoidMSpectral[{1, 1, 1}, 10]["M"] - 4/3

(* convergence in L, machine precision *)
TableForm[Table[With[{r = ellipsoidMSpectral[ax, L]},
{L, NumberForm[r["M"], 16], ScientificForm[Last[r["degreeNorms"]], 3],
r["minRho"]}], {L, 0, 24, 4}]]

(* 50 digits: compare two precisions, then print *)
ellipsoidMSpectral[ax, 48, WorkingPrecision -> 70, "ExtraNodes" -> 60]["M"] -
ellipsoidMSpectral[ax, 48, WorkingPrecision -> 60]["M"]
NumberForm[ellipsoidMSpectral[ax, 56, WorkingPrecision -> 60]["M"], 55]

(* the ellipsoid: support fraction about 0.75, upper below M *)
With[{r = ellipsoidMplusQP[ax, 4000]},
{r["n"], r["support"], NumberForm[r["lower"], 14], r["kktViolation"],
NumberForm[ellipsoidMplusUpper[r, ax]["upper"], 14]}]
\end{lstlisting}

\subsection*{Acknowledgement}
Anthropic Claude Opus 5.5 played a crucial role in turning 
my conjectures into theorems: to say that the AI model merely ``assisted'' 
would be an understatement.  It performed an absolutely essential research 
role throughout every aspect of the paper.  No AI funding whatsoever for 
this study was received.  My interest in maximal energy began years ago 
(\cite{Finch19}) and continues (\cite{Finch26}).


\end{document}